\documentclass[sn-jnl]{sn-jnl}% 
\usepackage[numbers]{natbib}
\usepackage{graphicx}% Para incluir gráficos
\usepackage{multirow}% Combinar filas en tablas.
\usepackage{amsmath,amssymb,amsfonts}% De la ams
\usepackage{amsthm}% Para teoremas ams.
\usepackage{mathrsfs}% Fuentes matem adicionales.
\usepackage[title]{appendix}% Para apéndices.
\usepackage{xcolor}% Para colores en el documento.
\usepackage{textcomp}% Para símbolos de texto
\usepackage{manyfoot}% Para pie de página y notas.
\usepackage{booktabs}% Para apariencia de tablas.
\usepackage{algorithm}% Para algoritmos.
\usepackage{algorithmicx}% Para algoritmos.
\usepackage{algpseudocode}% Paraalgoritmos.
\usepackage{listings}% Para fragmentos de código fuente.
\usepackage{mathtools}
\usepackage{csquotes} % Para gestionar citas
\DeclareGraphicsExtensions{.pdf, png} % Extensiones de gráficos
\usepackage{bbm} % Fuentes matemáticas especiales
\usepackage{enumerate} % Personalización de listas enumeradas
\usepackage{float} % Control de flotantes
\usepackage{cancel} % Tachado de texto en ecuaciones
\usepackage[font=bf, justification=centering]{caption} % Leyendas de figuras y tablas
\usepackage{geometry} % Personalización de márgenes
\usepackage{threeparttable} % Tablas con notas al pie
\usepackage{afterpage} % Ejecutar comandos después de la página actual
\usepackage{pdflscape} % Cambiar la orientación de la página a paisaje
\usepackage{placeins} % Controlar la posición de flotantes

\usepackage{titletoc} % Formato de contenidos
\usepackage{threeparttable}

\date{\today}

\newtheorem{theorem}{Theorem}%  meant for continuous numbers
\newtheorem{proposition}[theorem]{Proposition}% 
\newtheorem{definition}{Definition}%
\usepackage[numbers]{natbib}
\hypersetup{%
        hidelinks,
        breaklinks=true,
        plainpages=false,%
        citecolor=false,
        linkcolor=dalse,
        urlcolor=false,
        bookmarksopen=true,%
        bookmarksnumbered=false,%
        bookmarksdepth=5%
}

\date{October 30, 2023}
\begin{document}

\title[Article Title]{\textbf{
Time Scale Calculus: A New Approach to Multi-Dose Pharmacokinetic Modeling\textsuperscript{*}}\footnotetext{\textsuperscript{*}The authors did not receive support from any organization for the submitted work.}}

\author[1]{\fnm{Jos\'e Ricardo} \sur{Arteaga Bejarano}}\email{storresp@uchicago.edu}

\author*[2]{\fnm{Santiago} \sur{Torres}}\email{jarteaga@uniandes.edu.co}
\equalcont{These authors contributed equally to this work.}

\affil[1]{\orgdiv{Mathematics Department and Research Group in Mathematical and Computational Biology (BIOMAC)}, \orgname{Universidad de los Andes}, \orgaddress{\street{Cra 1 Nº 18A - 12}, \city{Bogot\'a}, \country{Colombia}}}

\affil*[2]{\orgdiv{Pre-Doctoral Fellow, The Pearson Institute, Harris School of Public Policy}, \orgname{University of Chicago}, \orgaddress{\street{1307 East 60th Street}, \city{Chicago}, \country{USA}}}


\abstract{In this paper, we use Time Scale Calculus (TSC) to formulate and solve pharmacokinetic models exploring multiple dose dynamics. TSC is a mathematical framework that allows the modeling of dynamical systems comprising continuous and discrete processes. This characteristic makes TSC particularly suited for multi-dose pharmacokinetic problems, which inherently feature a blend of continuous processes (such as absorption, metabolization, and elimination) and discrete events (drug intake).  We use this toolkit to derive analytical expressions for blood concentration trajectories under various multi-dose regimens across several flagship pharmacokinetic models. We demonstrate that this mathematical framework furnishes an alternative and simplified way to model and retrieve analytical solutions for multi-dose dynamics. For instance, it enables the study of blood concentration responses to arbitrary dose regimens and facilitates the characterization of the long-term behavior of the solutions, such as their steady state.   }



\maketitle

%%%%%%%%
%%%%%%%%
%%%%%%%%
\newpage
%% ====== SECTION 1 ===== %%
\section{Introduction}

Time Scale Calculus (TSC) is an emerging field in mathematics that aims to unify the theory of differential and difference equations. While differential and difference equations require modeling time as exclusively a continuous or discrete quantity, TSC considers dynamics happening on unrestricted time domains called time scales.  Hence, the core innovation of TSC lies in its capacity to simultaneously address continuous and discrete dynamics, offering an improved modeling framework for analyzing systems where these processes coexist. This approach was first introduced in the seminal works of \cite{Hilger1988} and \cite{HILGER_1997} and has attracted growing interest for its applicability across multiple disciplines.\\

This paper highlights the potential applications of TSC in Pharmacokinetic/Pharmacodynamic (PKPD) and Physiologically Based Pharmacokinetic (PBPK) modeling. In particular, we study blood concentration dynamics resulting from multiple administrations of orally and intravenously administered drugs.  Understanding these dynamics is crucial for medical practice, as drugs are typically administered in a series of repetitive doses, not in single doses. TSC proves invaluable in examining these dynamics, offering a versatile framework that naturally accommodates the dual nature of drug administration: the continuous processes of absorption, metabolization, and elimination, alongside the discrete events of drug intake. Consequently, TSC provides a coherent language for integrating these processes into models, enhancing their interpretability, and expanding their applications and complexity.\\

The primary application examined herein is modeling multiple-dose dynamics through the Bateman function \cite{niazi1979textbook}. The Bateman function, which articulates the dynamics following some orally administered single doses, emerges from the solution of a second-order linear ordinary differential equation or, equivalently, from a system of two first-order linear equations, that rationalize well-known physiological processes.  Thus, to extend this model to accommodate multiple doses, we employ TSC to transition from ordinary linear equations to dynamic equations - their generalization within the TSC framework
- enabling the incorporation of dosing timings into the dynamics. This methodology demonstrates how TSC facilitates the derivation of analytical solutions for established problems, such as the ``Equi-dosing'' regime - wherein a constant drug quantity is administered at uniform intervals - as studied originally studied in \cite{dost1953blutspiegel}, and further developed by \cite{hof2021exact}. However, we also show that TSC allows us to go far beyond this case, yielding solutions for arbitrary dosing regimes where the amount of drug administered and the interval between doses can vary arbitrarily. This generalization is important, so long equi-dosing regimens are ``a difficult, if not
impossible, feat to accomplish in any circumstances, be it a
clinical or preclinical study or, even more so, in the real case of a
patient taking a drug long term.''\cite[p.432]{brocks2010rate}.\\

We also use the multiple-dose dynamics of the Bateman function obtained via TSC to showcase results concerning drug accumulation. For instance, we establish mathematical properties inherent to these dynamics, such as asymptotic periodicity and the convergence to a ``Steady State''—a phase where blood concentration dynamics fluctuate within a predetermined range.  Additionally, we formulate a new mathematical definition of the ``Steady State'', and use it to provide explicit formulas for the asymptotic bounds of the dynamics.\\

Finally, we showcase the versatility of TSC in extending its application to various established PKPD and PBPK models. Consequently, we formulate multiple-dose dynamics via TSC of bolus injections through intravascular routes and models incorporating finite absorption times, as discussed in \cite{macheras2020revising,chryssafidis2022re}. Following the insights of
\cite{brocks2010rate}, our findings highlight the profound influence of administration routes and the underlying assumptions about drug absorption and elimination on drug accumulation kinetics in plasma.\\

Our paper is the first use of Time Scale Calculus in PKPD and PBPK applications. This novel application broadens TSC's scope beyond its established utility in areas such as population dynamics and biology (\textit{e.g.}, \cite{berger2013application,hu2012,spedding2003taming,thomas2009spray}), economics (\textit{e.g.}, 
\cite{Atici2006,tisdell2008basic,atici2009comparison,dryl2014calculus}
), physics (\textit{e.g.},  \cite{Ahlbrandt2000}), machine learning (\textit{v.g.}, \cite{Seiffertt2009}), and robotics (\textit{e.g.}, \cite{Akn2019ControlOW}). \\

The closest work to our proposal is that of  \cite{hof2021exact}, which presents exact solutions for equi-dosing regimens within multi-dose pharmacokinetic models incorporating transit compartments. Similarly, the research conducted by  \cite{savva2021} and \cite{savva2022} provides exact solutions for multiple intermittent intravenous infusions within a single-compartment model. Our research distinguishes itself from these precedents in three fundamental aspects: First, we show that the multi-dose dynamics can be naturally encapsulated within the TSC framework. Specifically, dynamic equations —a generalization of differential and difference equations proposed by TSC - allow us to integrate absorption, elimination, and dose intakes into a single initial value problem. Hence, this approach differs from traditional modeling practices where researchers typically amalgamate multiple conventional ordinary differential equations describing post-dose dynamics. Second, we show that the dynamic equation language to PKPD and PBPK allows a simpler formulation and solution to complex dosing regimens. We demonstrate this by deriving analytical solutions for arbitrary dosing schedules, thus improving upon the equi-dosing assumptions discussed in prior works. Third, our modeling framework is valid for models different from the standard transit compartment systems. An example is the multi-dose dynamics of Physiologically Based Finite Time  Pharmacokinetic (PBFTPK) models explored in this paper. With this application, we contribute to the PBFTPK literature  by providing multiple-dose dynamics  where only single-dose dynamics have been developed (see \cite{macheras2020revising,chryssafidis2022re}).\\

 In addition, we provide several new results concerning drug accumulation. First, we offer a new mathematical formulation of the ``steady state'' for multi-dose dynamics based on TSC. We also provide alternate proofs of several well-established drug accumulation results in medical practice, as described in \cite{brocks2010rate}. Finally, we offer a novel result that asserts that a physician can  - theoretically -  formulate dosage regimens for blood concentrations to stay within any desired range asymptotically.\\

The remainder of this paper is organized as follows.  Section \ref{S: Model} introduces Time Scale Calculus, derives multiple-dose dynamics stemming from the Bateman Function, and describes its properties. Section \ref{S: Applications} provides additional applications of TSC to other PKPD and PBPK models. Section \ref{S: Conclusion} concludes.

%% ====== SECTION 3 ===== %%
\section{An application of Time Scale Calculus: Exploring multi-dose dynamics of the Bateman Function} \label{S: Model}

In this section, we employ Time Scale Calculus (TSC) to deliver analytical solutions for multiple-dose dynamics resulting from the Bateman function under different dosing regimens. The Bateman function is a one-compartment model disposition with first-order absorption and elimination rate, frequently applied in pharmacokinetic analyses of oral dosing. As a solution to an ordinary differential equation initial value problem, we use the Bateman function to exemplify the value of TSC in transitioning from single-dose to multiple-dose dynamics via dynamic equations. \\

\subsection{The Bateman Function}

The Bateman Function describes the blood concentration dynamics following a single oral dose intake. As shown in  \cite{khanday2016}, the Bateman Function is the solution of a system of linear ordinary differential equations modeling drug absorption and circulation through the gastrointestinal tract. In fact, the Bateman Function is also a solution to a second-order linear differential equation, as shown in \eqref{EQ: FAN ORAL}.\\

Following \cite{FAN_2014}, the concentration of many orally-administered drugs in plasma over time after a single dose, denoted as $x(t)$, is well approximated in empirical data by the following \textquote{Bateman Function}:

\begin{equation}\label{EQ: FAN ORAL} % p.o. = Per os = Orally
x(t)=\dfrac{\kappa_a F d}{V\,(\kappa_a-\kappa_e)}\; \left(e^{-\kappa_e t}-e^{-\kappa_a t}\right)
\end{equation}

\noindent where $\kappa_a,\kappa_e$ represent the compound absorption and elimination rate constants, respectively, and $F,V,d$ account for absolute bioavailability, the volume of distribution, and the amount of the oral dose, respectively.  It is usually the case that  \(\kappa_{a} > \kappa_{e} > 0\), meaning the absorption of the drug into the central compartment is faster than the elimination process. However, in some cases  \(0<\kappa_{a} < \kappa_{e}\), which is known as the flip-flop situation \citep{MORTENSEN2008}. Although it is possible to deal with the case where both parameters are equal, we will not refer to this case in our theory since it is uncommon to be found in practice.\\

More broadly, and following equation \eqref{EQ: FAN ORAL}, the drug blood concentration dynamics resulting from the administration of a broad spectrum of medications is the solution to the following second-order initial value problem:

\begin{equation} \label{Eq:FAN IVP}
\begin{cases}
x''(t)+(\kappa_a + \kappa_e) x'(t) + \kappa_e\kappa_a x(t) = 0\\
x'(0) = {\kappa_a \,\gamma \,d }\\
x(0)=0\\
\end{cases}    
\end{equation}

\noindent where parameter $\gamma$ condenses different combinations of biological parameters depending on the route of administration\footnote{
For example, \cite{FAN_2014} indicates that for an oral administration (P.O.) of the drug $\gamma = F/V$ $[ml^{-1}]$ where $F$ is the bioavailability and $V$ is the volume of distribution, while for an intravenous administration (I.V.) of the drug, $\gamma$ involves the rate constants governing formation, $\kappa_{f}$, and elimination, $\kappa_{m}$, as well as  systemic availability, $F_{H}(m)$,  which is the ratio of the amount of metabolite leaving the liver to the amount of metabolite formed, and distribution volume $V_{m}$ for the metabolite.}.\\

Alternatively, this particular second-order linear equation can be transformed into a system of first-order linear equations by introducing an auxiliary function $y(t)$. The function \(y(t)\) can be biologically interpreted as representing the quantity of the drug present in organs responsible for absorption at a time $t$. This reformulation leads to the following dynamical system 

%----------------
\begin{equation}\label{Eq:DYN_SYS}
\begin{split}
\begin{cases}
y'(t) = -\kappa_{a} y(t)\\
x'(t) = \kappa_{a}\,\gamma \,y(t) - \kappa_{e} x(t)\\
y(0) = y_{0} = d;\\ 
x(0) = x_{0} = 0
\end{cases}
\end{split}
\end{equation}
 \noindent with $\kappa_{a} > 0$, $\kappa_{e} > 0$, and $\kappa_a \neq \kappa_e$.\\

From a mathematical perspective, we have a system of two autonomous linear first-order differential equations with constant coefficients, for which a closed-form solution is known. For one, the differential equation governing $y(t)$ represents the process by which the drug amount in the gastrointestinal tract is absorbed, wherein $\kappa_a$ is the absorption rate constant. For another, the dynamics of $x(t)$  encapsulate the aftermath of the absorption and elimination processes: $\kappa_{a}\gamma y(t)$ denotes the concentration of the drug instantaneously absorbed at time $t$ (inflow), while $- \kappa_{e} x(t)$ indicates the drug's elimination, wherein $\kappa_e$ is the elimination rate constant. The system is coupled with initial conditions $y(0) = y_{0} = d$, and $x(0) = x_{0} = 0$, where $d$ represents the quantity of the drug administered orally. This particular choice of initial conditions allows the solution to match equation \eqref{EQ: FAN ORAL}.\\

Transitioning from single-dose to multi-dose dynamics necessitates a framework capable of accommodating dose intake mechanics beyond what is typically offered by ordinary differential equations.  This requirement is met through the development of dynamic equations, which extend traditional differential equations to incorporate discrete events, such as drug intakes, into the model. In the subsequent section, we delineate the foundational elements of this theory, setting the stage for its application to the Bateman function.

\subsection{Time Scale Calculus Preliminaries}

Time Scale Calculus (TSC) is an emerging field in mathematics that aims to unify the theory of differential and difference equations by enabling the study of dynamics on more general time domains called time scales. As an introduction, we briefly overview the main definitions of TSC following \cite{Bohner2001} and \cite{Bohner2003}.\\

%-------------------- SUBSUBSECTION  1.1
\subsubsection{Basic Definitions}
As the name suggests, this theory aims to generalize discrete and continuous dynamical systems into more general time sets called time scales. Specifically, a time scale is defined as follows:\\

%%====== DEFINITION 1
\begin{definition} % DEF 1.1
A \emph{Time Scale} \(\mathbb{T}\) is an arbitrary nonempty closed subset of \(\mathbb{R}\).
\end{definition}
\vspace{2mm}
For instance, \(\mathbb{R},\; \mathbb{Z},\; \mathbb{N} \cup [0,1]\) are time scales.\\

%%====== DEFINITION 2
\begin{definition} % DEF 1.2
Let \(\mathbb{T}\) be a Time Scale. For \(t \in \mathbb{T}\), we define the \emph{forward jump operator} \(\sigma: \mathbb{T} \rightarrow \mathbb{T}\) by:
\[
\sigma(t):= \inf \{s \in \mathbb{T} : s>t\}
\]
Analogously, we define the \emph{backward jump operator}\\ \(\rho: \mathbb{T} \rightarrow \mathbb{T}\) by:
\[
\rho(t):= \sup \{s \in \mathbb{T} : s<t\}
\]
\end{definition}

\vspace{2mm}
Intuitively, the definition of time scales admits both discrete and continuous sets and sets with both components. As it is expected, dynamics will vary depending on the ``structure'' of the Time Scale. This leads to the following classification. \\

%%====== DEFINITION 3
\begin{definition} % DEF 1.3
Let \(\mathbb{T}\) a time scale. An element \(t \in \mathbb{T}\) is:
\begin{enumerate}%[label=(\roman*)]
\item
\emph{right-scattered} if \(t<\sigma(t)\).
\item 
\emph{right-dense} if \(t=\sigma(t)\).
\item 
\emph{left-scattered} if \(t>\rho(t)\).
\item 
\emph{left-dense} if \(t=\rho(t)\).
\item 
\emph{isolated} if \(\rho(t)<t<\sigma(t)\).
\item 
\emph{dense} if \(\rho(t)=t=\sigma(t)\)
\end{enumerate}
\end{definition}
\vspace{2mm}
In other words, the forward and backward operators determine whether, at a certain point, the dynamics resemble a discrete or a continuous scale. An alternative approach is to define another operator called the graininess function:
\vspace{2mm}
%%====== DEFINITION 4
\begin{definition} % DEF 1.4
Let \(\mathbb{T}\) a time scale. We define the \emph{graininess function} 
\(\mu : \mathbb{T} \longrightarrow \mathbb{R}_{\geq0}= [0, \infty)\) by:
$$\mu(t):= \sigma(t)-t$$
\end{definition}

Observe that if $\mu(t)=0$, then the time scale at $t$ is continuous, whereas if $\mu(t)>0$, then the time scale at $t$ behaves like a discrete set. Finally, we define a distinguished subset of every time scale that will be useful for other definitions.
\vspace{2mm}
%%====== DEFINITION 5
\begin{definition} % DEF 1.5
We define the subset \(\mathbb{T}^{\kappa}\) as
\[
\mathbb{T}^{\kappa}= \begin{cases}
\mathbb{T} \setminus (\rho(\sup(\mathbb{T})),\sup(\mathbb{T})] & \text{if } \sup(\mathbb{T})<\infty\\
\mathbb{T} & \text{if } \sup(\mathbb{T})=\infty\\
\end{cases}
\]
\end{definition}

%-------------SUBSUBSECTION 1.2
\subsubsection{Differentiation}
We introduce the extended notion of a derivative in an arbitrary time scale.\\

%%====== DEFINITION 6
\begin{definition} % DEF 1.6
Let \(f: \mathbb{T} \rightarrow \mathbb{R}\) be a function and \(t \in \mathbb{T}^{\kappa}\). We define \(f^{\Delta}(t)\) (provided it exists) to be the number that for every \(\epsilon>0\), there exists \(\delta>0\) such that there is a neighbourhood \(U=(t-\delta,t+\delta) \cap \mathbb{T}\) that satisfies:
\[
|(f(\sigma(t))-f(s))-f^{\Delta}(t)(\sigma(t)-s)| \leq \epsilon |\sigma(t)-s|, \quad \forall s \in U
\]
We call \(f^{\Delta}(t)\) the \emph{delta-derivative} (or Hilger Derivative) of \(f\) at \(t\).
\end{definition}
\vspace{2mm}
The following Theorem shows how this definition extends the notion of the regular derivative and discrete differences:\\

%%====== PROPOSITION 1
\begin{proposition} % THEO 1
Assume \(f:\mathbb{T} \rightarrow \mathbb{R}\) be a function and let \(t \in \mathbb{T}^{\kappa}\). Then the following are valid:

\begin{enumerate}[(i)]
\item 
if \(f\) is delta-differentiable at \(t\), it is continuous at \(t\).
\item 
if \(f\) is continuous and \(t\) is right-scattered, the \(f\) is delta-differentiable at \(t\) with: 
\[
f^{\Delta}(t)=\dfrac{f(\sigma(t))-f(t)}{\mu(t)}
\] 
\item 
if \(t\) is right-dense, then \(f\) is delta-differentiable at \(t\) if and only if the limit, 
\[
\lim \limits_{s \to t} \dfrac{f(t)-f(s)}{t-s}
\] 
exists as a finite number. In this case:
\[
f^{\Delta}(t)=\lim \limits_{s \to t} \dfrac{f(t)-f(s)}{t-s}=f'(t)
\]
\end{enumerate}
\end{proposition}

In summary, the Hilger derivative behaves like an ordinary derivative in dense domains and like a discrete difference in scattered domains. This enables a single and general  notion of “change” in time domains that exhibit both denseness and scatteredness properties.

% --------- SUBSECTION 3.1.3
\subsubsection{Dynamic Equations}

Like ordinary differential and difference equations, the TSC allows formulating equations that describe the dynamics of an object in an arbitrary time scale - these are called \textquote{dynamic equations}. A dynamic equation of order $p$ is a mathematical equation that relates a function and its delta derivatives up to order $p$ with respect to a single independent variable $t$. More precisely, it takes the form:

$$G(t,f(t),f^{\Delta}(t),\cdots, f^{\Delta^p}(t))=0$$

Additionally, these equations can be combined with initial and boundary conditions to produce specific solutions. The general solution for many dynamic equations is known and thoroughly addressed in \cite{Bohner2001, Bohner2003}.

%----------------------- SUBSECTION 3.2
\subsection{The equi-dosing regimen} \label{S: Equidose}

In this section, we apply Time Scale Calculus to model blood concentration levels resulting from multiple successive drug doses under a Bateman function. The key idea behind this proposal is that the specific dynamical system depicted in \eqref{Eq:DYN_SYS}, leading to the Bateman Function, results from assuming dynamics unfold on a quite specific time scale: $\mathbb{R}_{\geq0} = [0, \infty)$. However, these same dynamics can also potentially describe treatment responses in more general time scales. For instance, consider a simple scenario in which $d$ grams of a particular drug is administered every $\tau$ units of time, with the first dose being administered at time $t=0$. Thus, we can visualize the dynamics happening on the following time scale $\mathbb{T}(\tau)$:

\begin{equation}
\mathbb{T}(\tau)=\bigcup\limits_{n=1}^{\infty}[(n-1)\,\tau, n\, \tau] = \bigcup\limits_{n=1}^{\infty}I_{n}.  
\end{equation}

The intervals $I_{n} = [(n-1)\tau, n \tau]$, for $n=1,2,3,\dots$, represent the time lapses between each dose administration. The beginning of interval $I_n$ marks the moment the $n$-th dose is administered. In other words, the $n$-th dose is administered at the beginning of interval $I_n$, precisely at time $(n-1)\tau$. Because $d$ and $\tau$ are constant, this case is known in the literature as the ``Equi-dosing case'' \cite{hof2021exact}. \\

Leveraging TSC tools, we can embed the dynamical system introduced in \eqref{Eq:DYN_SYS} into a different time scale---$\mathbb{T}(\tau)$. Specifically, this generalization is possible via the Hilger derivative.  The exact formulation of this extension is as follows:

%% ============= %
\begin{equation}\label{Eq:NEW IVP}
\begin{split}
\begin{cases}    
y^{\triangle} (t) = -\kappa_{a}y(t)\\
x^{\triangle} (t) = \kappa_{a} \, \gamma \,y(t) -\kappa_{e}x(t)
\\
\textbf{with initial conditions:}
\\
y(0) = d\\
x(0) = 0 
\\
\textbf{and multiplicity conditions:}
\\
\lim \limits_{t \to n\tau^+} y(t)=\lim \limits_{t \to n\tau^-} y(t)+  \, d; \quad n=1, 2,\cdots, \infty\\
\lim \limits_{t \to n\tau^+} x(t)=\lim \limits_{t \to n\tau^-} x(t); \quad n=1, 2,\cdots, \infty
\end{cases}
\end{split}
\end{equation}

Equation \ref{Eq:NEW IVP} condenses the equi-dosing problem into a single initial value problem using dynamic equations. Note that when considering a single-dose scenario $(\tau \to \infty)$, the proposed model matches exactly \eqref{Eq:FAN IVP} since $\mathbb{T}(\infty)=[0,\infty)=\mathbb{R}_{\geq0}$, and thus, the Hilger derivative is the same as the ordinary derivative.\\

In addition to the initial conditions typically used in ordinary differential equations' initial value problems, multiple-dose dynamics require specifying conditions at dose intake times. Our model assumes that the concentration in the systemic circulation does not rise abruptly when a new dose is administered orally, hence the continuity condition 
 $$\lim \limits_{t \to n\tau^+} x(t)=\lim \limits_{t \to n\tau^-} x(t)$$

In terms of the model, this implies that the drug concentration at the time of a new intake is equal to the concentration observed after $\tau$ units of time have elapsed since the last dose administration. Similarly, we require that

$$\lim \limits_{t \to n\tau^+} y(t)=\lim \limits_{t \to n\tau^-} y(t)+  \, d$$

This condition can be interpreted as follows: when a new drug dose is administered, it immediately enters the gastrointestinal tract, causing an abrupt increase in the amount of medicine to be absorbed by that system by $d$ units. \\

The solution to the above system thus describes the dynamics of $x$ and $y$ after administering a dose $d$ every $\tau$ units of time. Since this model is derived from the \textquote{Bateman Function}, the resulting function for $x(t)$ can be described as its equi-multiple-dose version. The formal result is presented in the following theorem:\\

%%====== THEOREM 2
\begin{theorem}[The Equi-multiple-dose Bateman Function]\label{Th:multidose_solution}

Consider the Time Scale defined by $\mathbb{T}(\tau)=\bigcup\limits_{n=1}^{\infty} I_{n} = \bigcup\limits_{n=1}^{\infty} [(n-1)\tau, n \tau]$. Then the initial value problem presented in \eqref{Eq:NEW IVP} on Time Scale $\mathbb{T}(\tau)$ has a unique solution given by:

%% ================= %
\begin{equation}\label{Eq:SYSTEM_SOLUTIONS}
\begin{cases}
y(t)= \sum \limits_{n=1}^{\infty} \mathbbm{1}[t \in I_{n}] \overset{(n)}{y}(t); \quad \overset{(n)}{y}(t) =  d
\left(
\dfrac{1 - \alpha^{n}}{1 - \alpha}
\right)
e^{-\kappa_{a}(t - (n-1)\tau)}
\\
x(t)= \sum \limits_{n=1}^\infty \mathbbm{1}[t \in I_{n}] \overset{(n)}{x}(t); \quad \overset{(n)}{x}(t) =
C_{1} e^{-\kappa_{e}(t - (n-1)\tau)}
-
C_{2} e^{-\kappa_{a}(t - (n-1)\tau)}
\\
\text{where}
\\
C_{1} = \left(
\dfrac{\kappa_{a}\, \gamma \,d}{\kappa_{a} - \kappa_{e}}
\right)
\left(
\dfrac{1 - \beta^{n}}{1 - \beta}
\right)
\\[10pt] 
C_{2} = \left(
\dfrac{\kappa_{a}\, \gamma \, d}{\kappa_{a} - \kappa_{e}}
\right)
\left(
\dfrac{1 - \alpha^{n}}{1 - \alpha}
\right)
\\[10pt] 
\alpha = e^{-\kappa_{a} \tau}, 
\quad 0 < \alpha < 1
\\
\beta = e^{-\kappa_{e} \tau},
\quad 0 < \beta < 1  \\
\end{cases}    
\end{equation}

The function $x(t)$ is called the \textbf{Equi-multiple-dose Bateman Function}.
\end{theorem}
%% ============= %
\begin{proof}
See Appendix \ref{AP: generalized}.\\

\end{proof}

Theorem \ref{Th:multidose_solution} presents an explicit formula for the dynamics of multiple doses, assuming that both the dose amount and the interval between doses remain constant. The solutions are articulated as infinite sums of functions, where $\overset{\scriptstyle{(n)}}{x}(t)$ and $\overset{\scriptstyle{(n)}}{y}(t)$ denote the particular solutions within a cycle, that is, the period between the $n-th$ and $(n+1)-th$ dose.

\subsection{Arbitary dosing regimens}

In this section, we offer a direct extension of the model described in Section \ref{S: Equidose} when allowing for dosage regimens with uneven times between doses and possibly different grammages in the doses.\\

Let $t_n$ be the time the $n-$th dose is administered, and define $\tau_n = t_{n} - t_{n-1}$ be the time elapsed between doses. Let $d_n$ represent the amount of drug administered at time $t_n$. Then, the sequence  $\{(\tau_{n},d_n)\}_{n=1}^\infty$ characterizes any conceivable dosage plan.\\

We must consider the times scale generated by an arbitrary regimen to characterize dynamics in a general setting. Similar to the constant case, this can be formulated as

$$\mathbb{T}(\{\tau_n\}_{n=1}^\infty)=\bigcup \limits_{n=1}^\infty I_n(\tau_1,\cdots, \tau_n)=\bigcup \limits_{n=1}^\infty [t_{n-1}, t_{n}]  $$

Notice this is a generalization of $\mathbb{T}(\tau)$, which results from the special case when  $\tau_n=\tau$ for all $n$, so that $I_n(\tau_1,\cdots, \tau_n)=I_n=[(n-1)\tau,n\tau]$. \\

As with the equi-dosing plans, the corresponding dynamics are described by a system of dynamic equations. Concretely, for all $t \in \mathbb{T}(\{\tau_n\}_{n=1}^\infty)$, consider the multi-dose dynamics are given by the following system of dynamic equations,\\

%----
\vspace{2mm}
\begin{equation}\label{Eq:General IVP}
\begin{split}
\begin{cases}
y^{\triangle} (t) = -\kappa_{a}y(t)\\
x^{\triangle} (t) = \gamma \,\kappa_{a} \,y(t) -\kappa_{e}x(t)
\\
\textbf{with initial conditions:}
\\
y(0) = d_{1}\\
x(0) = 0 
\\
\textbf{and multiplicity conditions:}
\\ %1
\lim \limits_{t \to t_{n-1}^+} y(t) = \lim \limits_{t \to t_{n-1}^-} y(t)  +  \, d_{n}; \quad n=2, 3,\cdots, \infty
\\%2
\lim \limits_{t \to t_{n-1}^+} x(t) = \lim \limits_{t \to t_{n-1}^-} x(t); \quad n=2, 3,\cdots, \infty 
\end{cases}
\end{split}
\end{equation}
\vspace{2mm}

As in Section \ref{S: Equidose}, this dynamic system allows for a closed-form solution, which we present in Theorem \ref{Th:multidose_solution_general}\\

%%====== THEOREM 2
\begin{theorem}[The Generalized-multiple-dose Bateman Function]\label{Th:multidose_solution_general}

Provided that $\kappa_{a} \neq \kappa_{e}$ and given any arbitrary dosage schedule given by the sequence $\left\{(d_{n}, \tau_{n})\right\}_{n=1}^{\infty}$, the solution to the system of dynamic equations presented in \eqref{Eq:General IVP} is:

%--------
\begin{equation}
\begin{split}
y(t) &= \sum \limits_{n=1}^\infty \mathbbm{1}[t \in I_{n}] \overset{(n)}{y}(t)
\\
x(t) &= \sum \limits_{n=1}^\infty \mathbbm{1}[t \in I_{n}] \overset{(n)}{x}(t)
\end{split}
\end{equation}

where

\begin{equation}
\begin{split}
\stackrel{(n)}{y}(t) &=
\left(
\sum\limits_{i=1}^{n-1}\prod\limits_{j=i}^{n-1} d_{i}\alpha_{j} + d_{n}
\right)
e^{-\kappa_{a}(t - t_{n-1})}
\\
\stackrel{(n)}{x}(t) &=
C_{1}(n) \, e^{-\kappa_{e}(t - t_{n-1})}
-
C_{2}(n) \, e^{-\kappa_{a}(t - t_{n-1})}
\\
C_{1}(n) &= \dfrac{\kappa_{a}\, \gamma }{\kappa_{a} - \kappa_{e}}
\left(
\sum\limits_{i=1}^{n-1}\prod\limits_{j=i}^{n-1} d_{i}\alpha_{j} + d_{n}
\right)
+
\dfrac{\kappa_{a}\cdot \gamma}{\kappa_{a} - \kappa_{e}}
\left[
\sum\limits_{i=1}^{n-1}\prod\limits_{j=i}^{n-1} d_{i}\beta_{j}
- 
\sum\limits_{i=1}^{n-1}\prod\limits_{j=i}^{n-1} d_{i}\alpha_{j}
\right]
\\
C_{2}(n) &= \dfrac{\kappa_{a}\, \gamma\ }{\kappa_{a} - \kappa_{e}}
\left(
\sum\limits_{i=1}^{n-1}\prod\limits_{j=i}^{n-1} d_{i}\alpha_{j} + d_{n}
\right)
\\
\alpha_{s} &= e^{-\kappa_{a}\tau_{s}},\,
\beta_{s} = e^{-\kappa_{e}\tau_{s}},
\,
s=1,2,3,\dots, n\geq1.\\
\end{split}    
\end{equation}
\end{theorem}

The function $x(t)$ is called the \textbf{Generalized-multiple-dose Bateman Function}.

\begin{proof}
See Appendix \ref{AP: generalizedarbitary}.\\    
\end{proof}

\begin{figure}[H] % FIGURE 8
\centering
\caption{\label{fig:PERIODIC_VS_IRREGULAR_DOSE} Periodic Dose vs. Irregular Dose - Illustration}
\begin{threeparttable}
\includegraphics[scale=0.8]{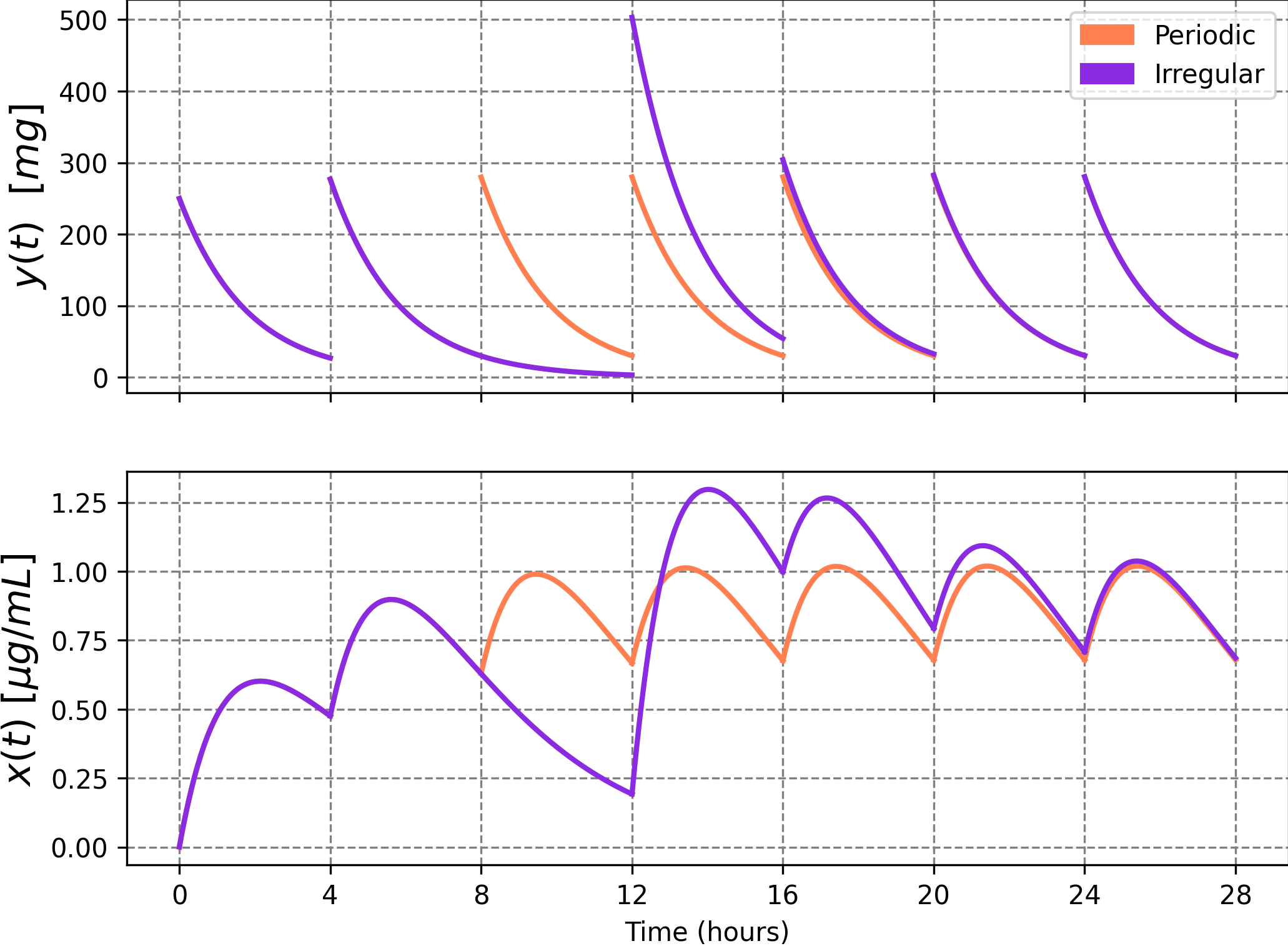}
\begin{tablenotes}
\item \textbf{Notes:} The top panel shows the amount of the drug in the intestinal tract, and the lower panel depicts the drug concentration in the bloodstream. 
\end{tablenotes}
\end{threeparttable}
\end{figure}

As illustrated in Figure \ref{fig:PERIODIC_VS_IRREGULAR_DOSE}, irregular dosing schedules can lead to dynamics that diverge significantly from the ``periodic'' solutions outlined in the previously discussed section. The blue line charts the blood concentration resulting from a consistent regimen of administering $250 \, mg$ of a hypothetical drug every 4 hours. Conversely, the orange line depicts a scenario where the patient skips the third dose and subsequently takes a double dose ($500 \, mg$) at the next scheduled time. While blood concentration levels eventually converge to a similar dynamic, this example underscores how deviations in dosing times and amounts can markedly affect the dynamics in the short term and may lead to prolonged variations if such patterns are consistently followed.

%% ====== FIGURE 8

%---------SUBSECTION 3.3
\subsection{Properties of the Multiple-dose Bateman Functions }

In this section, we exhibit properties of the Equi-multiple-dose Bateman Function. Similar properties can be established for the Generalized-multiple-dose Bateman Function. 

%------- SUBSUBSECTION 3.3.1
\subsubsection{Standard pharmacokinetic/pharmacodynamic quantities:}

We begin by calculating several typical pharmacokinetic/pharmacodynamic variables, such as the Area Under the Curve (AUC), maximum plasma concentration and the peak time in each cycle. Propositions \ref{Prop: AUC_In} through \ref{Prop: Maxcons} summarize these results. \\

%%====== PRPOSITION 3
\begin{proposition}\label{Prop: AUC_In}(Area under the curve - $\text{AUC}_{I_{n}}$)
The area under the concentration-time curve in the interval $I_{n} = [(n-1)\tau, n\,\tau]$, denoted as $\text{AUC}_{I_{n}}$, satisfies the following expression:
%---
\begin{equation}\label{Prop: Eq_AUC}
\text{AUC}_{I_{n}} = 
\dfrac{\kappa_{a}\, d\, \gamma}{\kappa_{a} - \kappa_{e}}
\left[
\dfrac{1}{\kappa_{e}}
\left(
1 - \beta^{n}
\right)
-
\dfrac{1}{\kappa_{a}}
\left(
1 - \alpha^{n}
\right)
\right]
\end{equation}
\end{proposition}
%---
\begin{proof}
See Appendix \ref{AP: Prop3}.\\
\end{proof}

%%====== PROPOSITION 4
\begin{proposition}\label{Prop: Bateman_AUC}
Let $\text{AUC}_{[0,\infty]}$ be the area under the concentration-time curve in the case of a single dose between $t=0$ and $t\to \infty$. The $\text{AUC}_{[0,\infty]}$ satisfies the following expression:
%---
\begin{equation}\label{Eq_AUC_mono-dosis}
\text{AUC}_{[0, \infty)} = 
\dfrac{\kappa_{a}\, d\, \gamma}{\kappa_{a} - \kappa_{e}}
\left[
\dfrac{1}{\kappa_{e}}
-
\dfrac{1}{\kappa_{a}}
\right]
\end{equation}
\end{proposition}

%----
\begin{proof}
See Appendix \ref{AP: Bateman_AUC}    
\end{proof}
\vspace{0.4cm}

%%====== PROPOSITION 5
\begin{proposition}\label{Prop: Maximum}
The time  at which the plasma concentration reaches its maximum value in period $I_{n} = [(n-1)\tau, n\,\tau]$  is denoted as $\stackrel{(n)}{t}_{\text{max}}$ and satisfies the equation:
\begin{equation}\label{Eq:t_max}
\stackrel{(n)}{t}_{\text{max}} =
(n-1)\tau +
\dfrac{1}{\kappa_{a} - \kappa_{e}}
\ln\left(
\dfrac{\kappa_{a} C_{2}}{\kappa_{e} C_{1}}
\right)
\end{equation}
which is always a positive number regardless of whether $\alpha > \beta$ or $\beta > \alpha$. Here, $C_1$ and $C_{2}$ are the constants defined in Theorem \ref{Th:multidose_solution}.
\end{proposition}
%---
\begin{proof}
See Appendix \ref{AP: Maximum}.\\
\end{proof}

\begin{figure}[H]% FIGURE 3
\caption{\label{fig: PK_dyn_1} The Equi-multiple-dose Bateman Function -Illustration.}
\begin{threeparttable}
\centering
\hspace*{-1cm}
\includegraphics[scale=0.7]{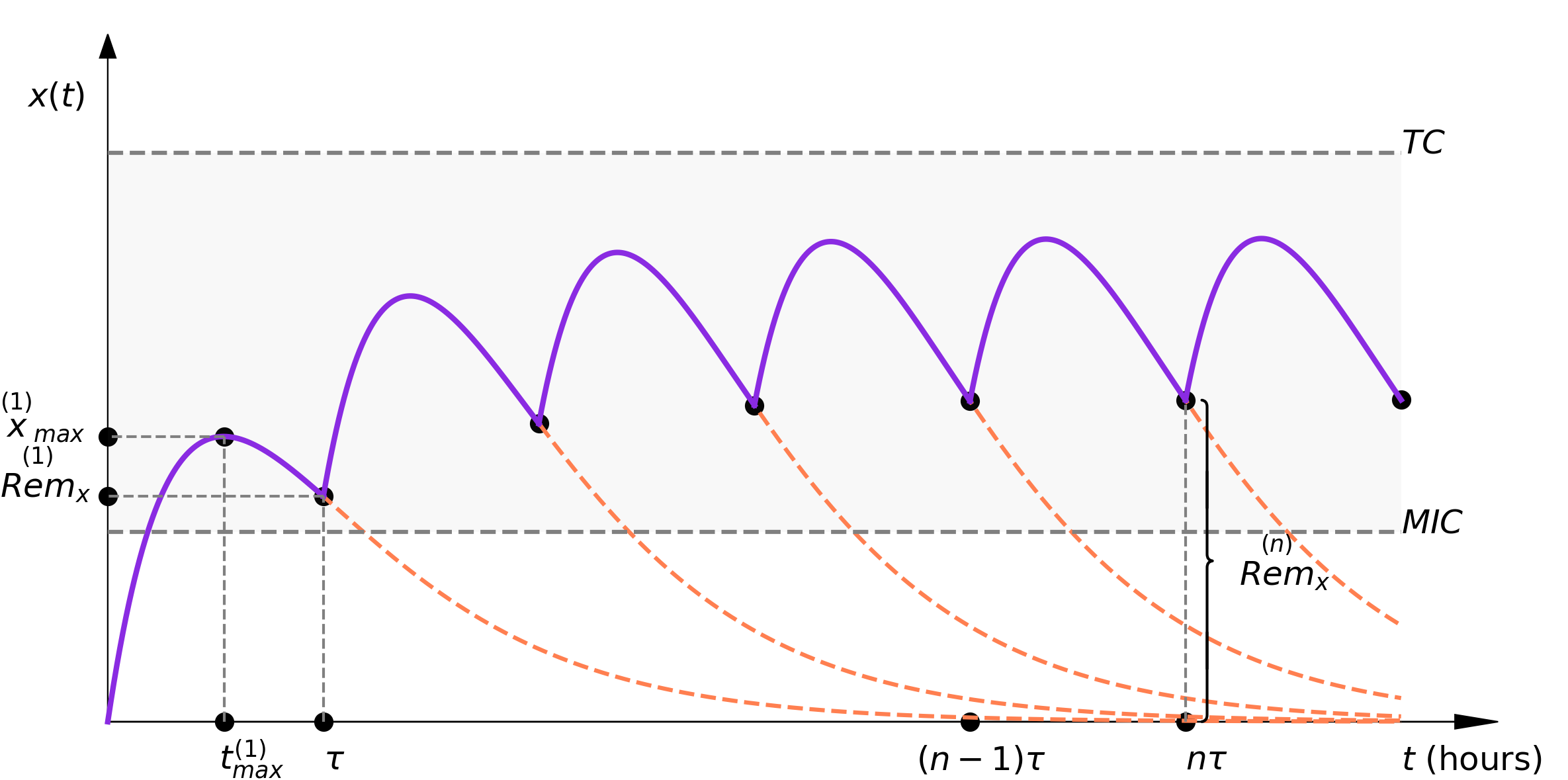} 
\end{threeparttable}
\end{figure}

%%====== PROPOSITION 6
\begin{proposition}\label{Prop: Maxcons}
The maximum plasma concentration in the period $I_{n} = [(n-1)\tau, n\tau]$, denoted as $\stackrel{(n)}{x}_{\text{max}}$, satisfies the following formula:
\begin{equation}
\stackrel{(n)}{x}_{\text{max}}
=
\stackrel{(n)}{x}\left(\stackrel{(n)}{t}_{\text{max}}\right)=
C_{1}
\left(
\frac{\kappa_{a}C_{2}}{\kappa_{e}C_{1}}
\right)^{-\dfrac{\kappa_{e}}{\kappa_a - \kappa_e}}
-
C_{2}
\left(
\frac{\kappa_{a}C_{2}}{\kappa_{e}C_{1}}
\right)^{-\dfrac{\kappa_{a}}{\kappa_a - \kappa_e}}
\end{equation}
where $C_{1}$ and $C_{2}$ are defined in Theorem \ref{Th:multidose_solution}.
\end{proposition}
%---
\vspace{2mm}
To summarize the previous results, we have graphically illustrated them in Figure \ref{fig: PK_dyn_1}.

%% ====== FIGURE 3

%----------- SUBSUBSECTION 3.3.2 
\subsubsection{Asymptotic behavior of the solutions.} \label{S: Asymptotic}

In this section, we investigate the plasma concentration after numerous dose administrations by observing the behavior of the sequence of functions $ \{\stackrel{(n)}{x}\}_{n =0}^\infty$ as $n \to \infty$. While these results pertain to the equi-dose regimen case, some can be extended to arbitrary dosage regimens in several situations.\\ 

%%====== THEOREM 7
\begin{theorem}[Asymptotic periodicity]\label{T: SSNew}

For a fixed choice of parameters ($\kappa_a,\kappa_e,\gamma$), $\kappa_a \neq \kappa_e$, the sequence of functions $\stackrel{(0)}{x}, \stackrel{(1)}{x}, \cdots, $ is asymptotically $\tau-$periodic, meaning

$$ \lim \limits_{n \to \infty}  \left[ \sup \limits_{t \in I_n} \; \Big |\stackrel{(n)}{x}(t)-\stackrel{(n-1)}{x}(t-\tau) \Big | \right]=0$$
\end{theorem}

\begin{proof}
See Appendix \ref{AP: SS}.\\    
\end{proof}

Theorem \ref{T: SSNew} asserts that the dynamics of plasma concentration stabilize into identical cycles after multiple doses have been administered. The predictability of concentration dynamics after several dose intakes suggests that solutions reach a \textquote{steady state}, where the blood concentration dynamics exhibit a regular pattern. A similar result can be established for many arbitrary dose regimens in Appendix \ref{A: arbitrary}. Consequently, we formalize the notion of a steady state as follows:\\

%%====== DEFINITION 7
\begin{definition}[Steady state] Let $\epsilon>0$ and let

$$N_\epsilon:= \min \Bigg \{ n \in \mathbb{N} \cup \{0\} : \forall m \geq n, \sup \limits_{t \in I_{m}} \; \Big |x^{(m)}(t)-x^{(m-1)}(t-\tau) \Big |<\epsilon \Bigg \}$$

We define the \textbf{$\epsilon$-steady-state} of the model as the time scale

$$\mathbb{T}_\epsilon^{ss}(\tau) := \bigcup_{n \geq N_\epsilon} I_n \subseteq \mathbb{T}(\tau)$$

Furthermore, we denote by 

$$\overset{(s.s.)}{x_{\epsilon}} := x  \Big \vert_{\mathbb{T}^{ss}_{\epsilon}(\tau)}$$

as the  \textbf{$\epsilon$-steady-state} dynamics.
\end{definition}

\vspace{2mm}

An immediate corollary of Theorem \ref{T: SSNew} is that \textquote{steady states} always exist for any conceivable dosage plan. More precisely, we know that a an $\epsilon$-steady-state always exist for any $\epsilon>0$, and for any set of parameters $\kappa_a,\kappa_e,\gamma$, provided that $\kappa_a \neq \kappa_e$.  Furthermore, this definition justifies studying the limiting behavior of standard pharmacokinetic/pharmacodynamic quantities by considering their limits when $n \to \infty$. \\

Additionally, asymptotic periodicity implies that it is possible to study the long-term behavior of patients following a specific dosage plan. In particular, it allows for examining the range of possible concentrations a patient may exhibit after being administered a particular drug dosage for a long period. \\

%%====== DEFINITION 8
\begin{definition}[Therapeutic Range]
Let $\epsilon>0$ and let

$$TR_\epsilon^{\text{s.s.}}=\text{Range}(\overset{(s.s.)}{x_\epsilon})$$ 

\noindent be the range of values a solution can exhibit when it has achieved its $\epsilon-$steady state. We define the therapeutic range or safety range as:
\begin{equation}
\text{TR}^{\text{s.s.}} 
=
[
\underline{\text{SS}},
\overline{\text{SS}}
]= \bigcap \limits_{\epsilon >0} TR_\epsilon^{\text{s.s.}}= \bigcap \limits_{\epsilon >0} \text{Range}(\overset{(s.s.)}{x_\epsilon})
\end{equation}
\end{definition}

\vspace{0.4cm}

In brief, the therapeutic or safety range is the set of possible concentrations an individual can exhibit when $t \to \infty$. We now find explicit formulas for $\underline{\text{SS}}$ and for $
\overline{\text{SS}}$. To achieve this, consider the following sequences indexed by $n$.\\

%%====== DEFINITION 9
\begin{definition}[Remainder]
We will refer to the plasma concentration of medication that remains in the circulatory system after the $n$-th period, $I_{n} = [(n-1)\tau, n \tau]$, as that period's remainder. We denote this quantity by $\stackrel{(n)}{\textbf{Rem}_{x}}$. Formally:  
%---
\begin{equation}\label{Eq:Remainder}
\stackrel{(n)}{\textbf{Rem}_{x}} := 
\stackrel{(n)}{x}(n\,\tau) =
\dfrac{\kappa_{a}\,d\,\gamma}{\kappa_{a} - \kappa_{e}}
\left[
\left(
\dfrac{1 - \beta^{n}}{1 -\beta}
\right)
\beta
-
\left(
\dfrac{1 - \alpha^{n}}{1 -\alpha}
\right)
\alpha
\right]    
\end{equation}
\end{definition}

\vspace{0.4cm}

Since both $0 < \alpha = e^{-\kappa_{a}\tau} < 1$ and $0 < \alpha = e^{-\kappa_{a}\tau} < 1$ as $n \to \infty$, the remainder converges to a quantity $\stackrel{(\infty)}{\textbf{Rem}_{x}}$. Moreover, by definition, this quantity constitutes the lower bound of the therapeutic range. Hence:  

\begin{equation}\label{Eq:Remainder_infinito}
\underline{\text{SS}}
=
\stackrel{(\infty)}{\textbf{Rem}_{x}} := 
\lim\limits_{n\to \infty} 
\stackrel{(n)}{\textbf{Rem}_{x}}
=
\dfrac{\kappa_{a}\,d\,\gamma}{\kappa_{a} - \kappa_{e}}
\left[
\left(
\dfrac{\beta }{1 -\beta}
\right)
-
\left(
\dfrac{\alpha}{1 -\alpha}
\right)
\right]    
\end{equation}

Furthermore, it can be verified that  $\stackrel{(n)}{\textbf{Rem}_{x}}$ is always a positive number for each $n$, regardless of whether $\alpha > \beta$ or $\beta > \alpha$. Likewise, it follows that   $\underline{\text{SS}}=\stackrel{(\infty)}{\textbf{Rem}_{x}} > 0$.\\

The following proposition carries out a similar procedure to find the upper bound of the therapeutic range\\

%%====== PROPOSITION 8
\begin{proposition}\label{Prop: Bounds}
The maximum plasma concentration in the steady state converges to a positive number and is given by the following expression:
\begin{equation}
\overline{\text{SS}} = \lim \limits_{n \to \infty}
\stackrel{(n)}{x}_{\text{max}} =
\frac{\kappa_{a}d\gamma}{\kappa_{a} - \kappa_{e}}
\left[
\frac{1}{1 - \beta}
\left(
\frac{\kappa_{a}(1-\beta)}{\kappa_{e}(1-\alpha)}
\right)^{-\frac{\kappa_{e}}{\kappa_{a} - \kappa_{e}}}
-
\frac{1}{1 - \alpha}
\left(
\frac{\kappa_{a}(1-\beta)}{\kappa_{e}(1-\alpha)}
\right)^{-\frac{\kappa_{a}}{\kappa_{a} - \kappa_{e}}}
\right]
\end{equation}
\end{proposition}
%---
\begin{proof}
See Appendix \ref{AP: Bounds}.\\
\end{proof}

Critically, plasma concentration asymptotic behavior depends on the dosage schedule $(d,\tau)$. Theorem \ref{Th:steady-state_bounds} demonstrates how different schedules result in different therapeutic ranges:\\

%%====== THEOREM 9
\begin{theorem}[Steady-state bounds]\label{Th:steady-state_bounds}
For a given vector physiological parameters, $(\kappa_{a}, \kappa_{e}, \gamma)$, let $\overline{SS}(d,\tau)$ and $\underline{SS}(d,\tau)$ be the upper and lower bounds of the therapeutic range in the steady-state seen as functions of $(d,\tau)$. \\
%---------
\begin{enumerate}
\item[(a)]
If the personalized dose $d$ is increased while maintaining the personalized dose regimen $\tau$ constant for a patient, then $\underline{SS} (d,\tau)$ increases and $\overline{SS} (d,\tau)$ increases.
\item[(b)]
If the time between doses $\tau$  is increased while keeping the dose $d$ constant, then $\underline{SS} (d,\tau)$ decreases with a limit of zero and $\overline{SS} (d,\tau)$ decreases up to a positive limit.
\end{enumerate}
\end{theorem}

\begin{proof}
See Appendix \ref{AP: generalized2}.\\    
\end{proof}

Moreover, we can show that the Therapeutic Range is well-defined in that $\overline{SS} > \underline{SS}$, as well as characterize how different dosage plans alter its width. Theorem \ref{Th:Therapeutic_range_width} exhibits these results\\

%%====== THEOREM 10
\begin{theorem}[Width of the Therapeutic Range]\label{Th:Therapeutic_range_width}
The maximum and minimum plasma concentrations, $\overline{SS}$ and $\underline{SS}$, in the steady-state always satisfy the following inequality:
\begin{equation}
\overline{SS}(d,\tau) > \underline{SS}(d,\tau)
\end{equation}
The width of the therapeutic range, defined as $\ell(d,\tau) := \overline{SS}(d,\tau) - \underline{SS}(d,\tau)$, increases with respect to the dose $d$ and the dosage regime $\tau$. Moreover, if the dose is kept constant ($d>0$), then,
\begin{equation}
\lim\limits_{\tau \to \infty}
\ell(d,\tau) =
\dfrac{\kappa_{a}\,d\, \gamma}{\kappa_{a} - \kappa_{e}}
\left[ 
\left(
\dfrac{\kappa_{a}}{\kappa_{e}}
\right)^{-\frac{\kappa_{e}}{\kappa_{a} - \kappa_{e}}}
-
\left(
\dfrac{\kappa_{a}}{\kappa_{e}}
\right)^{-\frac{\kappa_{a}}{\kappa_{a} - \kappa_{e}}}
\right]>0
\end{equation}
\end{theorem}

\begin{proof}
See Appendix \ref{AP: TWidth}. \\    
\end{proof}

Lastly, it is possible to characterize the area under the curve for the limiting cycle. This is given by

$$\text{AUC}^{s.s.}= \lim \limits_{n \to \infty} \text{AUC}_{I_{n}}$$

%%%
This quantity holds an interesting relationship with the single-dose AUC, as documented by \cite{brocks2010rate} and as formally proved in the following result.\\

%%====== THEOREM 11
\begin{theorem}[Equality of areas under the curves]\label{Th: equality_AUCs}
Let $\text{AUC}_{I_{n}}^{\,\,s.s.}$ be the area under the concentration-time curve in the period $I_{n} = [(n-1)\tau, n\tau]$. Then,
\begin{equation}
\text{AUC}_{[0,\infty)}
=
\text{AUC}^{\,\,s.s.}
\end{equation}
\end{theorem}

%---
\begin{proof}
See Appendix \ref{AP: equality_AUCs}. \\     
\end{proof}

The area under the concentration-time curve (AUC) is a critical measure quantifying medication absorption into the systemic circulation. Theorem \ref{Th: equality_AUCs} establishes that the AUC for a single dose is equivalent to the AUC for each cycle in a steady state. In other words, when prescribing an evenly spaced dosage plan with a fixed grammage, a medical practitioner can expect the AUC in the steady state, which is frequently unknown, to be equivalent to that of the single-dose scenario, which is frequently determined through clinical studies.\\

\begin{figure}[H]
\centering
\caption{\label{fig:FIG_4_MULTIPLE_DOSE} Asymptotic behavior of the Equi-multiple-dose Bateman Function}
\hspace*{-10mm}
\includegraphics[scale=0.7]{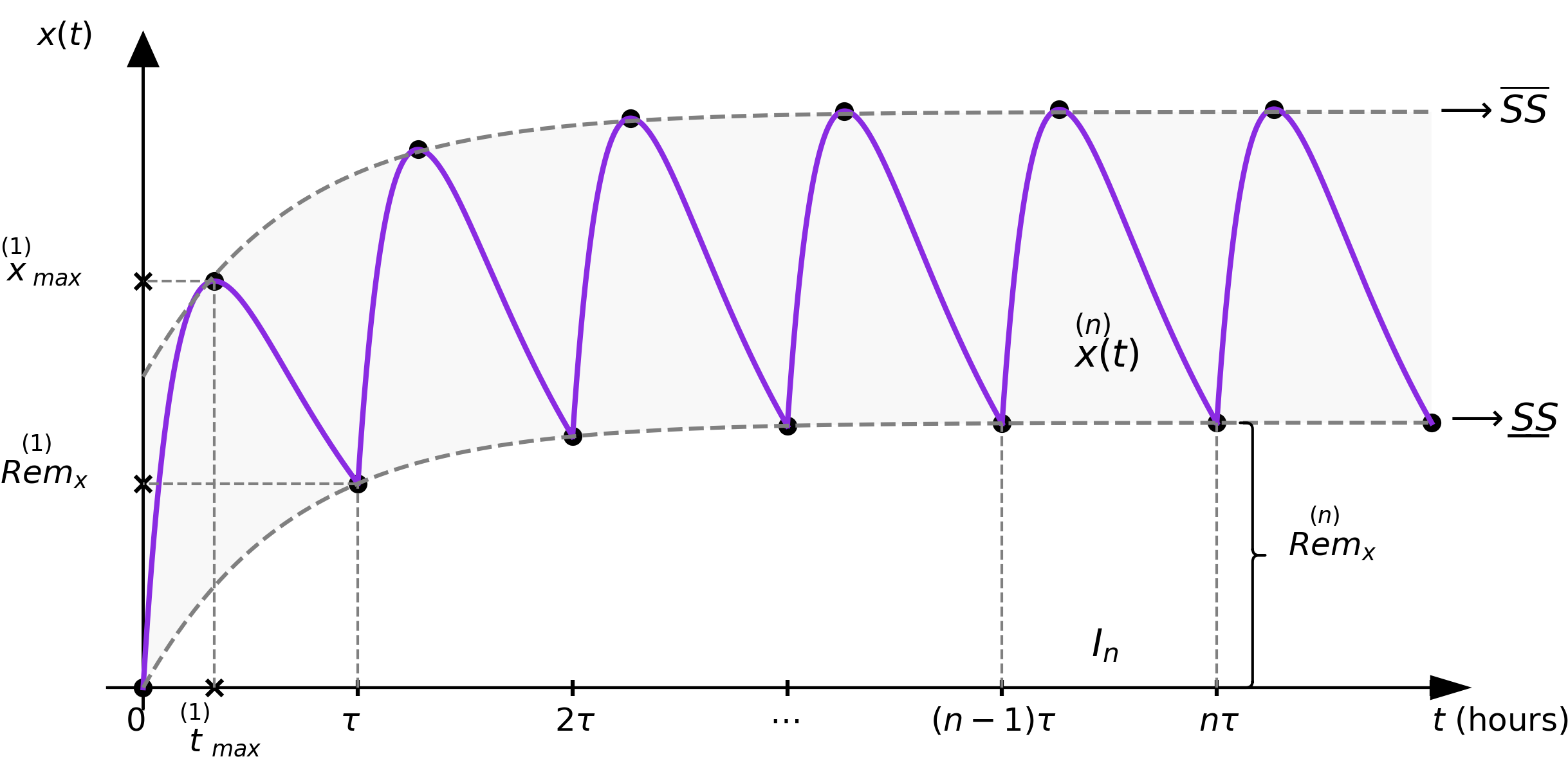}
\end{figure}

To summarize the previous results concerning the asymptotic behavior of the Generalized Function, we have graphically illustrated them in Figure \ref{fig:FIG_4_MULTIPLE_DOSE}.

%% ====== FIGURE 4

\subsubsection{Long-term treatments.}

Many long-term treatments require patients to adhere to a strict dosage regimen for several years, if not their entire lives. Some cases in question involve patients prescribed drugs for high blood pressure or HIV and non-disease treatments such as contraceptive pills. To facilitate the process for the patient, medical practitioners usually prescribe fixed-grammage doses to be taken at evenly spaced intervals. Furthermore, it is well known through practice that drug concentration in the blood stabilizes when prescribing these kinds of dosage schedules - a result we have now mathematically established to be always true in section \ref{S: Asymptotic}. Thus, the long-term behavior of these particular treatments can be well described using the asymptotic theory of the Equi-multiple-dose Bateman functions.\\

 In general, the objective of many of these treatments is to maintain concentration levels within a desired range $[\underline{R},\bar{R}]$, where $\bar{R}>\underline{R}>0$  \citep{jang2016, taddeo2020}). The lower bound is the Minimum Inhibitory Concentration (MIC), which represents the concentration above which the drug is effective. The upper bound is the Toxic Concentration (TC), which denotes concentrations above which the drug could potentially harm the user. As a result, an effective dosage plan maintains drug concentration levels after several doses within this range. In terms of the model, this means that

\begin{equation*}
\text{MIC}= \underline{R}
\leq\underline{\text{SS}}
< \overline{\text{SS}}
\leq \bar{R}= \text{TC}
\end{equation*}

From this description of the medical problem, two questions naturally arise: for a given set of parameters $(\kappa_a,\kappa_e,\gamma)$ characterizing a patient's kinetics, $i)$ does there always exist a dosage plan $(d,\tau)$ that results in a successful treatment for any  $[\underline{R},\bar{R}]$? and, $ii)$ provided there is a solution, how should a dosage plan be designed to be effective?\\

The answer to these questions is possible by studying the Equi-multiple-dose Bateman function's asymptotics. Specifically, we have established that the asymptotic behavior of the multi-dose dynamics is a function of $(d,\tau)$, meaning that the set of effective dosage schedules can be described as

\begin{equation*}
\mathcal{E}(\underline{R},\bar{R};\kappa_a,\kappa_e,\gamma)= \Big \{(d,\tau) \in \mathbb{R}_{>0} \times \mathbb{R}_{>0} : \underline{\text{SS}}(d,\tau;\kappa_a,\kappa_e,\gamma) \geq \underline{R}, \;\; \overline{\text{SS}}(d,\tau;\kappa_a,\kappa_e,\gamma) \leq \bar{R} \Big \}
\end{equation*}

\vspace{2mm}

%%====== THEOREM 12
\begin{theorem}\label{T: Existence}
For any $\bar{R}>\underline{R}>0$ and physiological parameters $(\kappa_a,\kappa_e,\gamma)$ with $\kappa_a \neq \kappa_e$, the set $\mathcal{E}(\underline{R},\bar{R};\kappa_a,\kappa_e,\gamma)$ is never empty.
\end{theorem}

\begin{proof}
See Appendix \ref{AP: Existence}.\\
\end{proof}

Theorem \ref{T: Existence} importance is twofold. From a medical standpoint, it assures us that there always exists a dosage plan—specifically, one featuring constant intervals between intakes and a fixed dosage—that ensures effective treatment for each patient. From a mathematical perspective, it guarantees that methods searching for these solutions can always find effective dosage plans for any patient and any imposed medical requirement.\\

In our particular framework, one could (theoretically) rely on our knowledge of the asymptotic behavior of the Equi-multiple-dose Bateman Equation to find effective dose schedules. For example, a health practitioner can fix desired long-run concentration levels $[\underline{\text{SS}}^*
,\overline{\text{SS}}^*] \subseteq [\underline{R},\bar{R}]$ for their patient. Assuming the biological parameters $(\kappa_a,\kappa_e,\gamma)$ were known, 
and effective dosage $(\tau^*,d^*)$ could be theoretically found by solving the  non-linear system of equations given by

\begin{equation*}
\begin{split}
\underline{\text{SS}}^*&=\underline{SS}(\tau^*,d^*,\kappa_a,\kappa_e,\gamma)\\
\overline{\text{SS}}^*&=\overline{SS}(\tau^*,d^*,\kappa_a,\kappa_e,\gamma)
\end{split}    
\end{equation*}

Unfortunately, given that such an approach is contingent upon the choice of a particular model and its unknown parameters, it does not offer a universally reliable tool in practice. However, the proof of Theorem \ref{T: Existence} does theoretically justify recurrent medical practices. For instance, it implies that the solutions to this problem can be (locally) expressed as functions of the parameters, namely:

\[
\begin{split}
\tau^*&=\tau^*(\underline{\text{SS}}^*,\overline{\text{SS}}^*,\kappa_a,\kappa_e,\gamma)\\
d^*&=d^*(\underline{\text{SS}}^*,\overline{\text{SS}}^*,\kappa_a,\kappa_e,\gamma)
\end{split} 
\]

In other words, the dependence of the solutions on the biological parameters endorses the well-known fact that not every drug dosage plan can be effective for every person.

%% ====== SECTION 4 ===== %%
\section{Other applications} \label{S: Applications}

This section showcases additional applications of Time Scale Calculus (TSC) to other PKPD and PBPK models. In particular, we explore multiple-dose dynamics of bolus injections through intravascular routes and models incorporating finite absorption times.

% ============ Example 1
\subsection{Additional application 1: Bolus injections via intravascular routes}

We consider a scenario where a drug is administered as a multidose bolus dose via an intravascular route. Let $x(t)$ denote the plasma concentration of the drug at time $t$, $t_n$ be the time the $n-$th dose is administered, $\tau_n = t_{n} - t_{n-1}$ be the time elapsed between doses, and $\delta_{n}$ denote the concentration of the drug delivered intravenously for dose $n$. In this context, the sequence  $\{(\tau_{n},\delta_n)\}_{n=1}^\infty$ outlines any potential dosing regimen. Therefore, the Time Scale suitable for this analysis is defined as:

$$\mathbb{T}(\{\tau_n\}_{n=1}^\infty)=\bigcup \limits_{n=1}^\infty I_n(\tau_1,\cdots, \tau_n)=\bigcup \limits_{n=1}^\infty I_n=\bigcup \limits_{n=1}^\infty [t_{n-1}, t_{n}]  $$

Given that the model only incorporates a central compartment to describe elimination dynamics, the dynamics of the model can be articulated through dynamic equations as follows:

\begin{equation}
\label{Eq:SEC4_EXAMPLE_1}
\begin{cases}
x^{\triangle}(t) = -\kappa_{e}x(t) \\
\textbf{with initial condition:}
\\
x(0) = \delta_{1}\\
\textbf{and multiplicity conditions:}
\\

\lim \limits_{t \to t_{n-1}^+} x(t)
=
\lim \limits_{t \to t_{n-1}^-} x(t)
+ \delta_{n}; 
\quad
n = 2, 3, \cdots , \infty 
\end{cases}
\end{equation}
where $\kappa_{e}$ is the elimination rate constant.\\

The analytical solution to the initial value problem formulated in \eqref{Eq:SEC4_EXAMPLE_1} is given in the following Theorem:\\

\begin{theorem}[The Generalized-multiple-dose bolus injection function]\label{Th:multidose_bolus}

Given any arbitrary dosage schedule given by the sequence $\left\{(\delta_{n}, \tau_{n})\right\}_{n=1}^{\infty}$, the solution to the system of dynamic equations presented in \eqref{Eq:SEC4_EXAMPLE_1} is:
%% ================= %
\begin{equation}\label{Eq:SEC4_EXAMPLE_1_SOL}
\begin{split}
x(t) &= \sum \limits_{n=1}^{\infty} \mathbbm{1}[t \in I_{n}] \overset{\scriptstyle{(n)}}{x}(t); 
\quad 
\overset{\scriptstyle{(n)}}{x}(t) =  
\left(
\sum\limits_{i=1}^{n-1}
\prod\limits_{j=i}^{n-1}\delta_{i}\beta_{j}
 + \delta_{n}
\right)
e^{-\kappa_{e}(t - t_{n-1})},
\quad
\beta_{s} = e^{-\kappa_{e}\tau_{s}} 
\end{split}
\end{equation}
\end{theorem}

\begin{proof}
See Appendix \ref{AP: BolusMulti}\\    
\end{proof}
% ---

\begin{figure}[H]
\centering
\caption{\label{fig:FIG_9_EXAMPLE_1} Multiple intravenous bolus dosing - Illustration} 
\begin{threeparttable}
\includegraphics[scale=0.1]{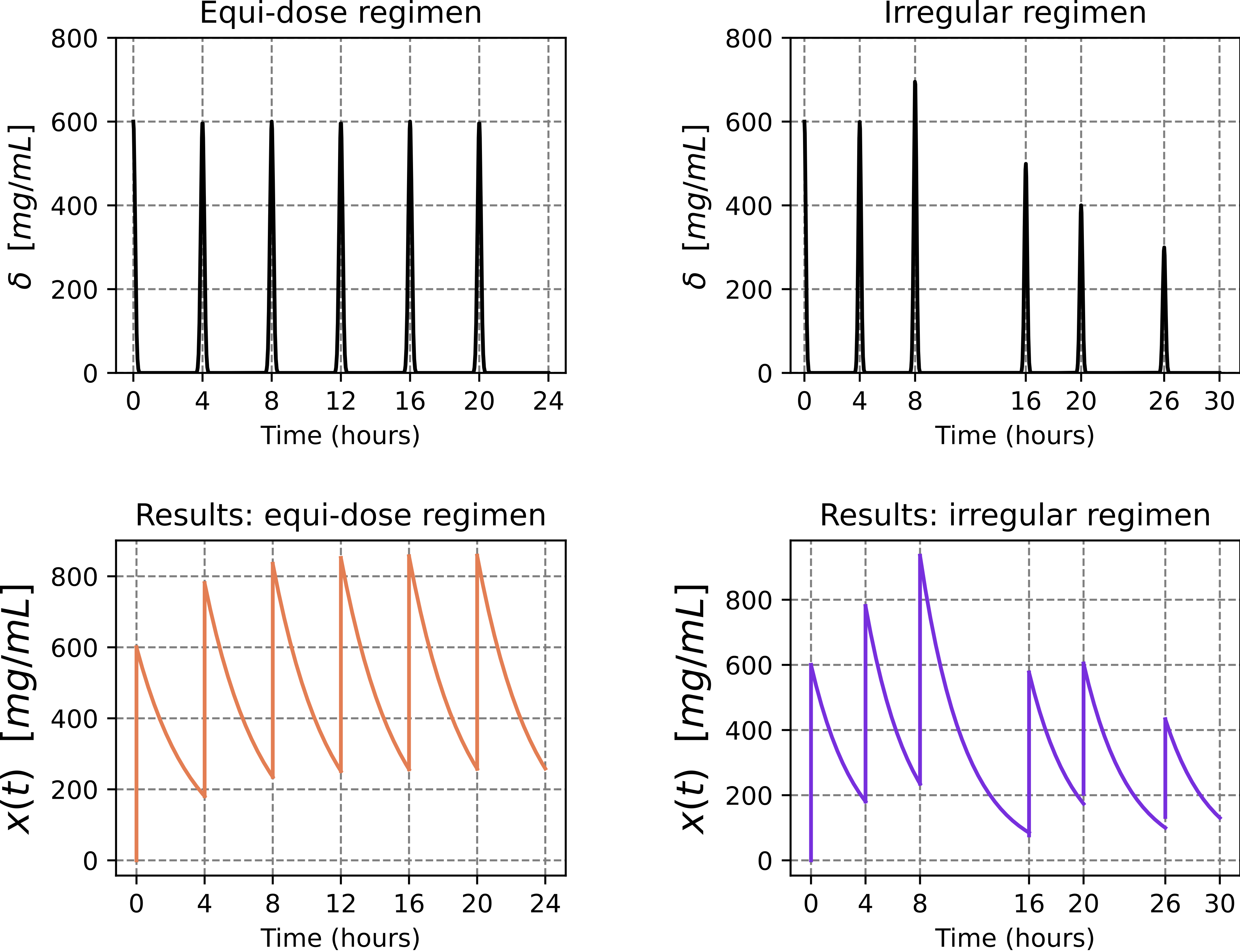}
\begin{tablenotes}
\item \textbf{Notes:}  The simulation shows an equi-dose regimen (left-hand side figures) with $d = 600 \, mg/mL$ and $\tau = 6$ hours, and an irregular regimen (right-hand side figures) with $d_{n} = [600, 600, 700, 500, 400, 300] \, mg/mL$ and $\tau_{n} = [4, 4, 8, 4, 6, 4]$ hours. In the simulation, $\kappa_{e} = 0.3838 \, h^{-1}$.
\end{tablenotes}
\end{threeparttable}
\end{figure}

Figure \ref{fig:FIG_9_EXAMPLE_1} uses Theorem \ref{Th:multidose_bolus} to demonstrate the dynamics of multiple doses under two distinct dosing regimens. The left panel displays the dynamics of traditional equi-dosing, while the right panel depicts an irregular dosing regimen. This comparison underscores how plasma concentration profiles can significantly differ due to variations in dosing schedules. As suggested by the functional form retrieved in Theorem \ref{Th:multidose_bolus}, this is a consequence of the dependency he dependency of the concentration during cycle $n$, $\overset{\scriptstyle{(n)}}{x}(t)$,  on the cumulative history of preceding doses $\left(\overset{\scriptstyle{(n)}}{x}(t)=\overset{\scriptstyle{(n)}}{x}(t;\tau_{n-1},\tau_{n-2},\cdots,\tau_1)\right)$.

\subsection{Application 2: Multple-dose Physiologically-Based Finite Time Pharmacokinetics models.}

Despite being a well-established model, the Bateman function has faced criticism, such as that articulated by \cite{macheras2019unphysical}, for its reliance on assumptions that may not accurately reflect physiological realities. One significant flaw identified by the authors is the presumption of an infinite absorption time for each dose administration. This does not align with the characteristics of some drugs known for their slow and incomplete absorption profiles, like those belonging to classes II,III, or IV. \\

Physiologically-based finite Time Pharmacokinetics (PBFTPK) models constitute a response to this critique, introducing the concept of a Finite Absorption Time (F.A.T). Essentially, the F.A.T. theory postulates two different phases following a dose intake. First, drug concentration dynamics are determined jointly by absorption and elimination forces, which can be referred to as the assimilation phase. However,  due to biological constraints, this phase is finite, and drug kinetics enter a second phase termed the clearance phase. In this stage, the drug concentration is governed solely by the elimination dynamics because absorption is no longer possible (\textit{e.g.}, the drug is no longer present in the gastrointestinal tract or beyond the respective absorptive sites). \\

Seminal works like \cite{macheras2020revising,chryssafidis2022re} have focused on single-dose PBFTPK models. We now extend these ideas to multiple dosing via TSC. Consider a regimen where a drug is administered in doses of $d_n$ milligrams every $\tau_n$ hours, with absorption occurring only within the first  $s_n$ hours, where $s_n \leq t_n$. As in other models, the sequence  $\{(\tau_{n},\delta_n)\}_{n=1}^\infty$  defines any conceivable dosing plan. However, the time scale here diverges from that in prior models, particularly in accommodating the intervals between doses as $ I_{n} = [t_{n-1}, s_{n}] \cup [s_{n}, t_{n}]$, with $[t_{n-1}, s_{n}]$ representing the assimilation phase and $ [s_{n}, t_{n}]$ the clearance phase. Thus, the new time scale incorporates both F.A.T. and dosing intervals, expressed as:

$$\mathbb{T}(\{s_n\}_{n=1}^\infty,\{\tau_n\}_{n=1}^\infty)=\bigcup \limits_{n=1}^\infty I_n=\bigcup \limits_{n=1}^\infty \left([t_{n-1}, s_{n}] \cup [s_{n}, t_{n}]\right)  $$

Within this framework, while maintaining the one-compartment model of Bateman with first-order absorption and elimination rates, the equations for multi-dose dynamics can be expressed as:

\begin{equation}\label{Eq:FAT_IVP}
\begin{split}
\begin{cases}    
% --
y^{\triangle} (t) =
-\kappa_{a}y(t)\\    
x^{\triangle} (t) =\kappa_{a} \, \gamma \,y(t) -\kappa_{e}x(t)
\\[20pt]
\textbf{with initial conditions:}
\\[10pt]
\overset{(1)}{y}(0) = d_{1}
\\
\overset{(1)}{x}(0) = 0
\\[10pt]
\textbf{and multiplicity conditions:}
\\[10pt]
y(t_{n-1}) = d_{n}; \quad n=1,2,\cdots,\infty  
\\
\lim \limits_{t \to t_{n-1}^+} x(t)
=
\lim \limits_{t \to t_{n-1}^-} x(t); \quad n=2, 3,\cdots,\infty    
\\
y(s_{n}) = 0; \quad n=1,2,\cdots,\infty  
\\
\lim \limits_{t \to s_{n}^+} x(t)
=
\lim \limits_{t \to s_{n}^-} x(t); \quad n=1,2,\cdots,\infty    
\end{cases}
\end{split}
\end{equation}

Equation \eqref{Eq:FAT_IVP} contains the initial value problem for the multiple-dose PBFTPK model. As in Section \ref{S: Equidose}, we impose multiplicity assumptions at dose intake times, namely, the first and second on the list. However, introducing F.A.T. requires the imposition of two additional multiplicity conditions. The third, and most critical, condition is expressed as:

$$y(s_{n}) = 0; \quad n=1,\cdots,\infty, $$

\noindent which implies that the absorption drug is depleted after time $s_n$.  Indeed, given that $y^{\triangle} (t) =
-\kappa_{a}y(t)=0$,  it follows that $y(t)=0$ for  $t \in [s_{n}, t_{n}]$, thus effectively encapsulating the theoretical dynamics of finite absorption.\\

The fourth condition:

$$\lim \limits_{t \to s_{n}^+} x(t)
=
\lim \limits_{t \to s_{n}^-} x(t); \quad n=1,\cdots,\infty   $$

ensures the continuity of blood concentration dynamics even after the drug has been fully absorbed.\\

The solution to this dynamic equation is detailed in Theorem \ref{Th:PBFTPK_2}:\\

\begin{theorem}[Multiple-dose PBFTPK drug concentration]\label{Th:PBFTPK_2}

Given any arbitrary dosage schedule given by the sequence $\left\{(\delta_{n}, \tau_{n})\right\}_{n=1}^{\infty}$,  and a sequence of F.A.T. $\{s_n\}_{n=1}^\infty$, the solution to the system of dynamic equations presented in \eqref{Eq:FAT_IVP} is:

%%%%%%%%%%%%%%%
\begin{equation}\label{Eq:FAT_IVP_SOL}
\begin{cases}
\, y(t)&= \sum \limits_{n=1}^{\infty} \mathbbm{1}[t \in I_{n}^{1}] \,\,\overset{\scriptstyle{(n)}}{y_{1}}(t), \quad
\overset{\scriptstyle{(n)}}{y_{1}}(t) =  d_{n}
e^{-\kappa_{a}(t - t_{n-1})} 
\\
x(t) &= \sum \limits_{n=1}^\infty \mathbbm{1}[t \in I_{n}^{1}] \,\,\overset{\scriptstyle{(n)}}{x_{1}}(t)+\mathbbm{1}[t \in I_{n}^{2}] \,\,\overset{\scriptstyle{(n)}}{x_{2}}(t); \quad 
\begin{cases}
\overset{\scriptstyle{(n)}}{x_{1}}(t) =
C_{1}(n) e^{-\kappa_{e}(t - t_{n-1})}
-
C_{2}(n) e^{-\kappa_{a}(t - t_{n-1})}
\\
\overset{\scriptstyle{(n)}}{x_{2}}(t) =
C_{3}(n) e^{-\kappa_{e}(t - s_{n})}  
\end{cases}
\\
\text{where}
\\[10pt]
I_{n}^{1} &= [t_{n-1}, s_{n}],\textbf{(Assimilation Phase)};
\quad 
I_{n}^{2} = [s_{n}, t_{n}],\textbf{(Clearance  Phase)}; \quad n = 1,\dots,\infty
\\[10pt]
C_{1}(n) &= 
\dfrac{\kappa_{a}\cdot\gamma}{\kappa_{a} - \kappa_{e}}
d_{n} +\dfrac{\beta}{B}\left[\dfrac{\kappa_{a}\cdot\gamma}{\kappa_{a} - \kappa_{e}}
(B - A)
\left(
\sum\limits_{i = 1}^{n-1}
\beta^{n-1-i}\,
d_{i}\right)
\right]
\\[10pt]
C_{2}(n) &= 
\dfrac{\kappa_{a}\cdot\gamma}{\kappa_{a} - \kappa_{e}} d_{n}
\\[15pt]
C_{3}(n) &= 
\dfrac{\kappa_{a}\cdot\gamma}{ \kappa_{a} - \kappa_{e}}
(B - A)
\left[
\sum\limits_{i = 1}^{n-1}
\beta^{n-1-i}\,
d_{i} 
\right]
\\[10pt]
\alpha &= e^{-\kappa_{a} \tau}, 
\quad 0 < \alpha < 1
\quad 
\beta = e^{-\kappa_{e} \tau},
\quad 0 < \beta < 1  
\\
A &= e^{-\kappa_{a} \sigma}, 
\quad 0 < A < 1
\quad 
B = e^{-\kappa_{e} \sigma},
\quad 0 < B < 1  
\\
\alpha^{n} &= e^{-\kappa_{a}t_{n}},
\quad
\beta^{n} = e^{-\kappa_{e}t_{n}},
\quad
\alpha^{n-1} A = e^{-\kappa_{a}s_{n}},
\quad
\beta^{n-1} B = e^{-\kappa_{e}s_{n}}; \quad n = 1,2,\dots,\infty
\\
\alpha^{0} &= \beta^{0}  = 1, \; t_0=0.    

\end{cases}    
\end{equation}
%%%%%%%%%%%%%%%
\end{theorem}

\begin{proof}
See Appendix \ref{AP: FATMulti}.\\ 

\end{proof}

Figure \ref{fig:FIG_11_EXAMPLE_3} showcases the dynamics of multiple doses incorporating Finite Absorption Time (F.A.T). The left panel depicts an equi-dosing scenario, where both the administered dose and the interval between dosing and the F.A.T remain constant. Conversely, the right panel introduces variations in these parameters across doses. A notable feature introduced by incorporating F.A.T is a distinctive kink in the blood concentration dynamics due to the solution provided in Theorem \ref{Eq:FAT_IVP_SOL} being a piecewise function\footnote{These characteristics appear in other works exploring blood concentration dynamics after single or multiple successive zero-order input rates under F.A.T schemes (see \cite{Alimpertis2023}). }. This change corresponds to the phase shift, in which blood concentration dynamics transit from being determined by both absorption and elimination forces, to being governed only by elimination.

\begin{figure}[H]
\centering
\caption{\label{fig:FIG_11_EXAMPLE_3} Multiple-dose F.A.T dynamics - Illustration} 
\begin{threeparttable}
\includegraphics[scale=0.1]{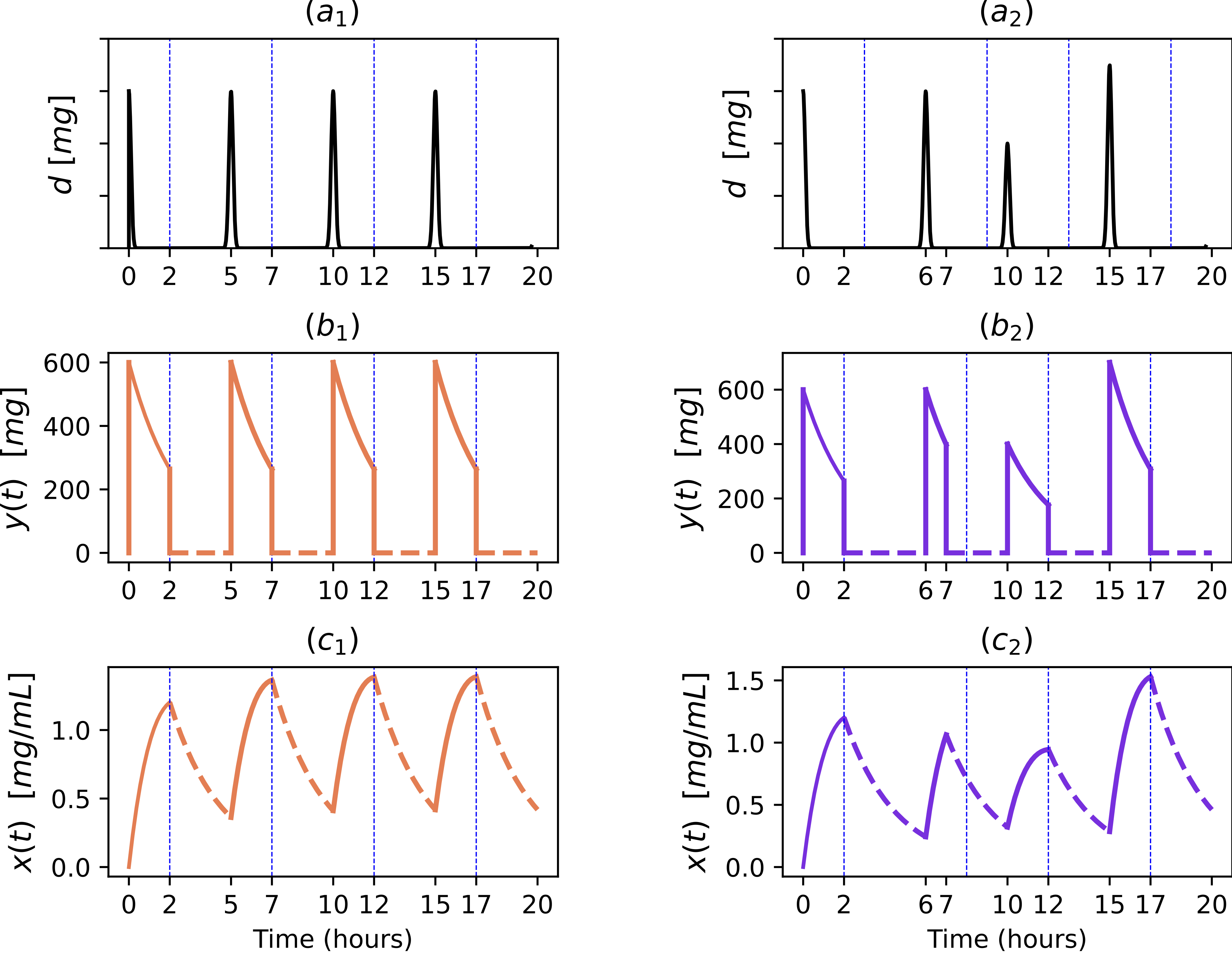}
\begin{tablenotes}
\item \textbf{Notes:}  
The simulation shows an equi-dose regimen (left-hand side figures) with $d = 600 \, mg$, $\tau = 5$ hours, $s_n = 2$ hours and an irregular regimen (right-hand side figures) with $d_{n} = [600, 600, 400, 700] \, mg$ and $\tau_{n} = [6, 4, 5, 5]$ hours, and $s_n = 2$ hours. In the simulation,
$\kappa_{a} = 0.42\, h^{-1}$
$\kappa_{e} = 0.4\, h^{-1}$,
$\gamma = 0.00449\, ml^{-1}$ . Vertical dotted lines indicate the Finite Absorption times. 
\end{tablenotes}
\end{threeparttable}
\end{figure}

%% ====== SECTION 6 ===== %%
\section{Conclusions} \label{S: Conclusion}

This paper applies Time Scale Calculus (TSC) to explore the dynamics of multiple drug doses, emphasizing how TSC's tools simplify our formulation of these complex dynamics. TSC offers a comprehensive framework ideally suited to the dual nature of drug administration, which encompasses the continuous processes of absorption, metabolism, and elimination, as well as the discrete instances of drug intake. Through applying TSC, we develop the multi-dose dynamics for several key pharmacokinetic (PK) models and establish analytical expressions for the trajectory of blood concentrations. Specifically, we formulate multiple-dose dynamics of several flagship PK models using TSC and derive new analytical formulas of the blood concentration dynamics. A critical advantage of TSC is that it enables the formulation of entire multi-dose dynamics as a simple, unique initial value problem. This simplification in the language allows us to handle intricate dosing patterns, like those featuring irregular dose timings and quantities, and to accommodate non-conventional absorption-elimination dynamics, such as those with a finite absorption period.\\

Although our discussion has been limited to a few simple examples that facilitate the exposition of the technique, TSC constitutes an ideal toolkit to formulate and solve more complex PKPD and PBPK models. For instance, by using the same techniques shown in this paper, it is possible to find analytical solutions for higher-order transit compartment models under unrestricted dosing regimens. Likewise, TSC offers possibilities to model coupled dynamical systems that incorporate continuous processes and discrete events. For example, it enables the simultaneous modeling of multi-dose drug administration aimed at bacterial eradication alongside models of bacterial population dynamics. Furthermore, TSC opens avenues for exploring optimization strategies for multiple-dose prescriptions via dynamic programming in time scales, a field well developed in recent years.

\newpage
\bibliographystyle{unsrt} % apalike, apacite, plainnat

\bibliography{sn-TSC}

\newpage
%%%%%%%%
%%%%%%%%
%%%%%%%%

\appendix

\begin{center}
    {\LARGE\bfseries Appendix for:} \\
    \vspace{0.4cm}
    
    {\LARGE\bfseries Time Scale Calculus: A New Approach to Multi-Dose Pharmacokinetic Modeling.}
\end{center}

\renewcommand{\partname}{}
\renewcommand{\thepart}{}
\setcounter{page}{0}
\thispagestyle{empty}

\part{} % Start the appendix part

\startcontents[appendix]
\printcontents[appendix]{}{1}{\noindent\textbf{\Large Contents}\vspace{0.5cm}\par}

\newpage

%% ====== SECTION 7 ===== %%
\section{Appendix - Asymptotic periodicity for arbitrary dose schedules.} \label{A: arbitrary}

Similarly to the equi-dose model, we can affirm the existence of asymptotic periodicity for non-equi-dose regimes, thereby establishing the concept of a steady state. Specifically, we can extend Theorem \ref{T: SSNew} through the following theorem.\\

\begin{theorem}[Asymptotic periodicity]\label{T: SSNewGeneral}

Consider an arbitrary dosage schedule $\{(d_n,\tau_n)\}_{n=1}^\infty$ such that $\lim \limits_{n \to \infty} d_n =d$ and $\lim \limits_{n \to \infty} \tau_n =\tau$, for some $d>0$ and $\tau>0$.\\

For a fixed choice of parameters ($\kappa_a,\kappa_e,\gamma$), $\kappa_a \neq \kappa_e$, the sequence of functions $\stackrel{(0)}{x}, \stackrel{(1)}{x}, \cdots, $ resulting from the Generalized Bateman function with schedule  $\{(d_n,\tau_n)\}_{n=1}^\infty$ is asymptotically $\tau-$periodic, meaning

$$ \lim \limits_{n \to \infty}  \left[ \sup \limits_{t \in I_n} \; \Big |\stackrel{(n)}{x}(t)-\stackrel{(n-1)}{x}(t-\tau) \Big | \right]=0$$

\end{theorem}

\begin{proof}

We begin this proof with a preliminary proposition:

\begin{proposition}
Consider the remainder sequences $\stackrel{(n)}{\textbf{Rem}_{x}}$ and $\stackrel{(n)}{\textbf{Rem}_{y}}$ presented in Theorem \ref{Teo:Formulas_recursivas}. Then

$$\limsup \limits_{n \to \infty} \; \stackrel{(n)}{\textbf{Rem}_{x}}= \overline{\textbf{Rem}}_x<\infty \quad \text{ and } \quad \limsup \limits_{n \to \infty}  \; \stackrel{(n)}{\textbf{Rem}_{y}}=\overline{\textbf{Rem}}_y<\infty $$

\end{proposition}

\begin{proof}

Since $\{\tau_n\}_{n=1}^\infty$ and $\{d_n\}_{n=1}^\infty$ are convergent sequences, they are also bounded. Let $\bar{\tau}=\max_n \tau_n$ and $\bar{d}=\max_n d_n$.\\

To prove the assertion, it suffices to show that both sequences are bounded under the given hypothesis. First, notice that for all $n$,  $\alpha_n=e^{-\kappa_a  \tau_n }$ is bounded. Let $\bar{\alpha}$ be this bound, with $\bar{\alpha} \in (0,1)$. This assertion is also true for $\beta_n=e^{-\kappa_e  \tau_n }$, whose bound we will denote $\bar{\beta} \in (0,1)$. \\

It follows that

\begin{align*}
\stackrel{(n)}{\textbf{Rem}_{y}} &= 
\sum\limits_{i=1}^{n}\prod\limits_{j=i}^{n} d_{i}\alpha_{j}
\leq \bar{d} \sum\limits_{i=1}^{n} \left( \bar{\alpha} \right)^{n-i}
= \bar{d} \sum\limits_{i=0}^{n-1} \left( \bar{\alpha} \right)^{i}
\leq \bar{d} \sum\limits_{i=0}^{\infty} \left( \bar{\alpha} \right)^{i}
=\dfrac{\bar{d}}{1-\bar{\alpha}}<\infty
\end{align*}

given the geometric series is convergent insofar as $\bar{\alpha} \in (0,1)$.\\

Similarly,

\begin{align*}
\stackrel{(n)}{\textbf{Rem}_{x}} &= 
\dfrac{\kappa_{a}\cdot \gamma}{\kappa_{a} - \kappa_{e}}
\left[
\sum\limits_{i=1}^{n}\prod\limits_{j=i}^{n} d_{i}\beta_{j}
- 
\sum\limits_{i=1}^{n}\prod\limits_{j=i}^{n} d_{i}\alpha_{j}
\right]\\
& \leq \dfrac{\kappa_{a}\cdot \gamma}{\kappa_{a} - \kappa_{e}} \left[
\sum\limits_{i=1}^{n}\prod\limits_{j=i}^{n} d_{i}\beta_{j}
+ 
\sum\limits_{i=1}^{n}\prod\limits_{j=i}^{n} d_{i}\alpha_{j}
\right]\\
& \leq \dfrac{\kappa_{a}\cdot \gamma \bar{d}}{\kappa_{a} - \kappa_{e}} \left[
\sum\limits_{i=1}^{n} \left( \bar{\beta} \right)^{n-i}
+ 
\sum\limits_{i=1}^{n} \left( \bar{\alpha} \right)^{n-i}
\right]\\
&=\dfrac{\kappa_{a}\cdot \gamma \bar{d}}{\kappa_{a} - \kappa_{e}} \left[
\sum\limits_{i=0}^{n-1} \left( \bar{\beta} \right)^{i}
+ 
\sum\limits_{i=0}^{n-1} \left( \bar{\alpha} \right)^{i}
\right]\\
&\leq \dfrac{\kappa_{a}\cdot \gamma \bar{d}}{\kappa_{a} - \kappa_{e}} \left[
\sum\limits_{i=0}^{\infty} \left( \bar{\beta} \right)^{i}
+ 
\sum\limits_{i=0}^{\infty} \left( \bar{\alpha} \right)^{i}
\right]\\
&\leq \dfrac{\kappa_{a}\cdot \gamma \bar{d}}{\kappa_{a} - \kappa_{e}} \left[
\dfrac{1}{1-\bar{\beta}}
+ 
\dfrac{1}{1-\bar{\alpha}}
\right]<\infty\\
\end{align*}

\end{proof}

Returning to the original proof, since $t_{n-1}=\tau_{n-1}+t_{n-2}$, it follows that

\footnotesize
\begin{align*}
\left| \stackrel{(n)}{x}(t)-\stackrel{(n-1)}{x}(t-\tau)\right| &= \Bigg \vert    \left[\left(\dfrac{\kappa_{a}\cdot \gamma \cdot(\stackrel{(n-1)}{\mathbf{Rem}_{y}} + d_{n})}{\kappa_{a} - \kappa_{e}} + \stackrel{(n - 1)}{\mathbf{Rem}_{x}}\right)e^{-\kappa_{e}(t - t_{n-1})}- 
\left(\dfrac{\kappa_{a}\cdot \gamma \cdot(\stackrel{(n-1)}{\mathbf{Rem}_{y}} + d_{n})}{\kappa_{a} - \kappa_{e}} \right)e^{-\kappa_{a}(t - t_{n-1})} \right]\\
&-\left[\left(\dfrac{\kappa_{a}\cdot \gamma \cdot(\stackrel{(n-2)}{\mathbf{Rem}_{y}} + d_{n-1})}{\kappa_{a} - \kappa_{e}} + \stackrel{(n - 2)}{\mathbf{Rem}_{x}}\right)e^{-\kappa_{e}(t-\tau - t_{n-2})}- 
\left(\dfrac{\kappa_{a}\cdot \gamma \cdot(\stackrel{(n-2)}{\mathbf{Rem}_{y}} + d_{n})}{\kappa_{a} - \kappa_{e}} \right)e^{-\kappa_{a}(t-\tau - t_{n-2})} \right] \Bigg \vert   \\
&= |A_n(t)\, e^{-\kappa_e t_{n-2}}-B_n(t)\, e^{-\kappa_a t_{n-2}}|\\
&\leq |A_n(t)|+|B_n(t)|
\end{align*}
\normalsize
    
where

\small
\begin{align*}
|A_n(t)|&= \left|  \left(\dfrac{\kappa_{a}\cdot \gamma \cdot(\stackrel{(n-1)}{\mathbf{Rem}_{y}} + d_{n})}{\kappa_{a} - \kappa_{e}} + \stackrel{(n - 1)}{\mathbf{Rem}_{x}}\right)e^{-\kappa_{e}(t - \tau_{n-1})}- \left(\dfrac{\kappa_{a}\cdot \gamma \cdot(\stackrel{(n-2)}{\mathbf{Rem}_{y}} + d_{n-1})}{\kappa_{a} - \kappa_{e}} + \stackrel{(n - 2)}{\mathbf{Rem}_{x}}\right)e^{-\kappa_{e}(t-\tau)}\right|\\
&\leq \left|  \left(\dfrac{\kappa_{a}\cdot \gamma \cdot(\stackrel{(n-1)}{\mathbf{Rem}_{y}} + d_{n})}{\kappa_{a} - \kappa_{e}} + \stackrel{(n - 1)}{\mathbf{Rem}_{x}}\right) - \left(\dfrac{\kappa_{a}\cdot \gamma \cdot(\stackrel{(n-2)}{\mathbf{Rem}_{y}} + d_{n-1})}{\kappa_{a} - \kappa_{e}} + \stackrel{(n - 2)}{\mathbf{Rem}_{x}}\right) \right|+\left|e^{-\kappa_{e}(t - \tau_{n-1})}- e^{-\kappa_{e}(t-\tau)}\right|\\
&= \left|  \left(\dfrac{\kappa_{a}\cdot \gamma \cdot(\stackrel{(n-1)}{\mathbf{Rem}_{y}} + d_{n})}{\kappa_{a} - \kappa_{e}} + \stackrel{(n - 1)}{\mathbf{Rem}_{x}}\right) - \left(\dfrac{\kappa_{a}\cdot \gamma \cdot(\stackrel{(n-2)}{\mathbf{Rem}_{y}} + d_{n-1})}{\kappa_{a} - \kappa_{e}} + \stackrel{(n - 2)}{\mathbf{Rem}_{x}}\right) \right|+\left|e^{\kappa_{e}\tau_{n-1}}- e^{\kappa_{e}\tau}\right| e^{-\kappa_{e}t}
\end{align*}
\normalsize

and

\small
\begin{align*}
|B_n(t)|&= \left|\left(\dfrac{\kappa_{a}\cdot \gamma \cdot(\stackrel{(n-1)}{\mathbf{Rem}_{y}} + d_{n})}{\kappa_{a} - \kappa_{e}} \right)e^{-\kappa_{a}(t - \tau_{n-1})}  -\left(\dfrac{\kappa_{a}\cdot \gamma \cdot(\stackrel{(n-1)}{\mathbf{Rem}_{y}} + d_{n})}{\kappa_{a} - \kappa_{e}} \right)e^{-\kappa_{a}(t - \tau)} \right|\\
&\leq \left|\left(\dfrac{\kappa_{a}\cdot \gamma \cdot(\stackrel{(n-1)}{\mathbf{Rem}_{y}} + d_{n})}{\kappa_{a} - \kappa_{e}} \right)  -\left(\dfrac{\kappa_{a}\cdot \gamma \cdot(\stackrel{(n-1)}{\mathbf{Rem}_{y}} + d_{n})}{\kappa_{a} - \kappa_{e}} \right) \right|+ \left|e^{-\kappa_{a}(t - \tau_{n-1})}-e^{-\kappa_{a}(t - \tau)} \right|\\
&= \left|\left(\dfrac{\kappa_{a}\cdot \gamma \cdot(\stackrel{(n-1)}{\mathbf{Rem}_{y}} + d_{n})}{\kappa_{a} - \kappa_{e}} \right)  -\left(\dfrac{\kappa_{a}\cdot \gamma \cdot(\stackrel{(n-1)}{\mathbf{Rem}_{y}} + d_{n})}{\kappa_{a} - \kappa_{e}} \right) \right|+ \left|e^{\kappa_{e}\tau_{n-1}}- e^{\kappa_{e}\tau}\right| e^{-\kappa_{e}t}
\end{align*}
\normalsize

So that

\begin{align*}
\sup \limits_{t \in I_n} \; |A_n (t)| &\leq  \left|  \left(\dfrac{\kappa_{a}\cdot \gamma \cdot(\stackrel{(n-1)}{\mathbf{Rem}_{y}} + d_{n})}{\kappa_{a} - \kappa_{e}} + \stackrel{(n - 1)}{\mathbf{Rem}_{x}}\right) - \left(\dfrac{\kappa_{a}\cdot \gamma \cdot(\stackrel{(n-2)}{\mathbf{Rem}_{y}} + d_{n-1})}{\kappa_{a} - \kappa_{e}} + \stackrel{(n - 2)}{\mathbf{Rem}_{x}}\right) \right|+\left|e^{\kappa_{e}\tau_{n-1}}- e^{\kappa_{e}\tau}\right|\\
\sup \limits_{t \in I_n} \; |B_n (t)| &\leq \left|\left(\dfrac{\kappa_{a}\cdot \gamma \cdot(\stackrel{(n-1)}{\mathbf{Rem}_{y}} + d_{n})}{\kappa_{a} - \kappa_{e}} \right)  -\left(\dfrac{\kappa_{a}\cdot \gamma \cdot(\stackrel{(n-1)}{\mathbf{Rem}_{y}} + d_{n})}{\kappa_{a} - \kappa_{e}} \right) \right|+ \left|e^{\kappa_{e}\tau_{n-1}}- e^{\kappa_{e}\tau}\right|
\end{align*}

and, consequently

$$\sup \limits_{t \in I_n}  \left| \stackrel{(n)}{x}(t)-\stackrel{(n-1)}{x}(t-\tau)\right| \leq \sup \limits_{t \in I_n} \; |A_n (t)| + \sup \limits_{t \in I_n} \; |B_n (t)| =C_n
 $$

Now note that

\begin{align*}
 \limsup \limits_{n \to \infty} \; \sup \limits_{t \in I_n} \;  |A_n (t)| &\leq   \left|  \left(\dfrac{\kappa_{a}\cdot \gamma (\overline{\textbf{Rem}}_y + d)}{\kappa_{a} - \kappa_{e}} + \overline{\textbf{Rem}}_x\right) - \left(\dfrac{\kappa_{a}\cdot \gamma (\overline{\textbf{Rem}}_y + d)}{\kappa_{a} - \kappa_{e}} + \overline{\textbf{Rem}}_x\right) \right|+\left|e^{\kappa_{e}\tau}- e^{\kappa_{e}\tau}\right|=0\\
  \limsup \limits_{n \to \infty} \; \sup \limits_{t \in I_n} \;  |B_n (t)| &\leq \left|\left(\dfrac{\kappa_{a}\cdot \gamma \cdot(\overline{\textbf{Rem}}_y + d)}{\kappa_{a} - \kappa_{e}} \right)  -\left(\dfrac{\kappa_{a}\cdot \gamma \cdot\overline{\textbf{Rem}}_y + d)}{\kappa_{a} - \kappa_{e}} \right) \right|+\left|e^{\kappa_{e}\tau}- e^{\kappa_{e}\tau}\right|=0
\end{align*}

Hence 

$$ \limsup \limits_{n \to \infty} \; C_n =0$$

Invoking the squeeze theorem, it follows that

$$ \lim \limits_{n \to \infty}  \left[ \sup \limits_{t \in I_n} \; \Big |\stackrel{(n)}{x}(t)-\stackrel{(n-1)}{x}(t-\tau) \Big | \right]=0$$
\end{proof}

\newpage
%% ====== SECTION 8 ===== %%
\section{Appendix - Omitted proofs - Bateman function}

\subsection{Proof of Theorem \ref{Th:multidose_solution}.} \label{AP: generalized}

\begin{proof}
We will use mathematical induction.
\begin{itemize}
\item
\textbf{$n=1$}. In this case, the Hilger derivative inside interval $I_{1} = [0, \tau] $ coincides with the usual derivative. Hence, the solution of the first equation of the system \eqref{Eq:NEW IVP} under the first initial condition  is,
\begin{equation}
\begin{split}
\dfrac{dy}{dt} &= -\kappa_{a} y \implies 
y(t) =  d\,  e^{-\kappa_{a} t}
\end{split}    
\end{equation}
Using this solution we can find the solution of the second equation of the system \eqref{Eq:NEW IVP} under the second initial condition is,
\begin{equation}
\begin{split}
\dfrac{dx}{dt} &= \kappa_{a}\, \gamma \, y(t) - \kappa_{e} x(t)\\
\implies &
x'(t) + \kappa_{e} \, x(t) = \kappa_{a}\, \gamma\, d \,e^{-\kappa_{a} t}
\end{split}    
\end{equation}
This is a non-homogeneous first-order linear differential equation with constant coefficients. Its integrating factor is $\mu = e^{\kappa_{e} \, t}$. Multiplying for the integrating factor and integrating with respect to $t$ we obtain,
\begin{equation}
\begin{split}
e^{\kappa_{e}\, t} \, x(t) &= 
\kappa_{a}\, \gamma \, d \, 
\dfrac{e^{-(\kappa_{a} - \kappa_{e}) t}}{-(\kappa_{a} - \kappa_{e})} + C
\end{split}    
\end{equation}
Now using the initial conditions: $y(0) = \gamma \, d$, $x(0) = 0$ we can find the constant of integration,
\begin{equation}
C = \dfrac{\kappa_{a}\, \gamma \, d}{\kappa_{a} - \kappa_{e}}    
\end{equation}
So, the solution of the second equation is,
\begin{equation}
x(t) = 
\dfrac{\kappa_{a}\, \gamma \, d}{\kappa_{a} - \kappa_{e}}
\left(
e^{-\kappa_{e} t}
-
e^{-\kappa_{a} t}
\right)
\quad 
t \in [0, \tau]
\end{equation}
\item 
\textbf{$n>1$.} We proceed by mathematical induction. \\

\textbf{Induction hypothesis.}\\
\\
For $t \in I_{n} = [(n-1)\tau, n\, \tau]$ and $n > 1 $ equations \eqref{Eq:SYSTEM_SOLUTIONS} are satisfied.\\

%%%%%%%%%%%
\textbf{Induction thesis (what must be proved)}
\\
For $t \in I_{n+1} = [n\tau, (n+1)\, \tau]$ the initial value problem,
%---
\begin{equation}\label{Eq:SYSTEM_INDUCTION_THESIS}
\begin{cases}
\stackrel{(n+1)\triangle}{y} (t) = -\kappa_{a} \, \stackrel{(n+1)}{y}(t)\\
\stackrel{(n+1)\triangle}{x} (t) = \kappa_{a}\, \gamma \, \stackrel{(n+1)}{y} (t) -\kappa_{e} \, \stackrel{(n+1)}{x} (t)\\
\text{with I.C.} 
\\
\overset{\scriptstyle{(n+1)}}{y}(n\tau) = \overset{\scriptstyle{(n)}}{y}(n \tau) +  d
\\
\overset{\scriptstyle{(n+1)}}{x}(n\tau) = \overset{\scriptstyle{(n)}}{x}(n \tau)
\end{cases}
\end{equation}

has explicit solution:

\begin{equation}
\begin{cases}
\overset{(n+1)}{y}(t) = \gamma \, d
\left(
\dfrac{1 - \alpha^{n+1}}{1 - \alpha}
\right)
e^{-\kappa_{a}(t - n\tau)}
\\
\overset{(n)}{x}(t) =
\tilde{C}_{1} e^{-\kappa_{e}(t - n\tau)}
-
\tilde{C}_{2} e^{-\kappa_{a}(t - n\tau)}
\\
\text{where}
\\
\tilde{C}_{1} = \left(
\dfrac{\kappa_{a}\, \gamma \,d}{\kappa_{a} - \kappa_{e}}
\right)
\left(
\dfrac{1 - \beta^{n+1}}{1 - \beta}
\right)
\\[10pt] 
\tilde{C}_{2} = \left(
\dfrac{\kappa_{a}\, \gamma \, d}{\kappa_{a} - \kappa_{e}}
\right)
\left(
\dfrac{1 - \alpha^{n+1}}{1 - \alpha}
\right)
\end{cases}    
\end{equation}
To prove the solution for $\stackrel{(n+1)}{y}(t)$ holds, we will start by finding an equivalent expression for the first initial condition. For the induction hypothesis, we have,
\begin{equation}
\begin{split}
\stackrel{(n+1)}{y}(n\tau) &= \stackrel{(n)}{y}(n\tau) + \gamma \, d
\\
&= \gamma \, d
\left(
\dfrac{1-\alpha^{n}}{1-\alpha} 
\right)e^{-\kappa_{a}\tau}
+ \gamma\, d
\\
&= \gamma \, d
\left(
\dfrac{1-\alpha^{n}}{1-\alpha} 
\right)\alpha
+ \gamma\, d
\\
&= \gamma \, d
\left[
\left(
\dfrac{1-\alpha^{n}}{1-\alpha} 
\right)\alpha
+ 1
\right]
\\
\stackrel{(n+1)}{y}(n\tau) 
&=
\gamma \, d
\left(
\dfrac{1-\alpha^{n+1}}{1-\alpha} 
\right)
\end{split}    
\end{equation}
Since the Hilger derivative in the interior of the interval $I_{n+1}$ is the derivative in the usual sense we solve the initial value problem,
\begin{equation}
\begin{split}
\dfrac{d}{dt}\stackrel{(n+1)}{y}(t) &= -\kappa_{a} \stackrel{(n+1)}{y}(t)
\\
\implies
\stackrel{(n+1)}{y}(t)
&=
\stackrel{(n+1)}{y}(n\tau) e^{-\kappa_{a}(t - n\tau)}
\\
\stackrel{(n+1)}{y}(t)
&=
\gamma \, d
\left(
\dfrac{1-\alpha^{n+1}}{1-\alpha} 
\right)e^{-\kappa_{a}(t - n\tau)}
\end{split}    
\end{equation}
Now, we will prove that the solution for $\stackrel{(n+1)}{x}(t)$ is satisfied in the interval $I_{n+1}$.
Using the induction hypothesis again and also taking into account the structure of the differential equation which is: non-homogeneous, first-order linear with constant coefficients and the non-homogeneous part is a multiple of the previously found solution for $\stackrel{ (n+1)}{y}(t)$ (exponential form), we can suppose that the solution has the form
\begin{equation}
\stackrel{(n+1)}{x}(t)
=
\tilde{C}_{1}
e^{-\kappa_{e}(t-n\tau)}
-
\tilde{C}_{2}
e^{-\kappa_{a}(t-n\tau)}
\end{equation}
where we only have to find the value of the constants, $\tilde{C}_{1}, \tilde{C}_{1}$. For this we will use the initial conditions of the problem and the induction hypothesis.
\begin{equation}\label{Eq:primera_C1_C2}
\begin{split}
\stackrel{(n+1)}{x}(n\tau) &=
\stackrel{(n)}{x}(n\tau)
\\
\implies
\tilde{C}_1 - \tilde{C}_2
&=
C_{1}e^{-\kappa_{e}\tau}
-
C_{2}e^{-\kappa_{a}\tau}
\\
\tilde{C}_1 - \tilde{C}_2
&=
C_{1}\beta
-
C_{2}\alpha
\end{split}    
\end{equation}
On the other hand, deriving $\stackrel{(n+1)}{x}(t)$ with respect to $t$ and evaluating in $n\tau$ we obtain a second equation for the constants we want to find.
\begin{equation}\label{Eq:segunda_C1_C2}
\begin{split}
\left.\dfrac{d}{dt}\stackrel{(n+1)}{x}(t)\right|_{t=n\tau}  
&=
-\kappa_{e}\tilde{C}_{1}
+
\kappa_{a}\tilde{C}_{2}
=
\kappa_{a}  \stackrel{(n+1)}{y}(n\tau) - \kappa_{e}\stackrel{(n+1)}{x}(n\tau)
\\
-\kappa_{e}\tilde{C}_{1}
+
\kappa_{a}\tilde{C}_{2}
&=
\kappa_{a} \stackrel{(n+1)}{y}(n\tau) - \kappa_{e}\stackrel{(n)}{x}(n\tau)
\\
-\kappa_{e}\tilde{C}_{1}
+
\kappa_{a}\tilde{C}_{2}
&=
\kappa_{a}
\left(
\dfrac{1-\alpha^{n+1}}{1-\alpha}
\right)
- \kappa_{e}\stackrel{(n)}{x}(n\tau)
\\
-\kappa_{e}\tilde{C}_{1}
+
\kappa_{a}\tilde{C}_{2}
&=
\kappa_{a}
\left(
\dfrac{1-\alpha^{n+1}}{1-\alpha}
\right)
- \kappa_{e}
\left(
C_{1}\beta - C_{2}\alpha
\right)
\end{split}    
\end{equation}
Now, multiplying the first equation, \eqref{Eq:primera_C1_C2}, on both sides by $\kappa_{e}$ and adding term by term with the second equation, \eqref{Eq:segunda_C1_C2}, and simplifying we have,
\begin{equation}
\tilde{C}_{2} = 
\dfrac{\kappa_{a}\gamma\, d}{\kappa_{a} - \kappa_{e}}
\left(
\dfrac{1-\alpha^{n+1}}{1 - \alpha}
\right)
\end{equation}
Finally, substituting this result in the first equation \eqref{Eq:primera_C1_C2} and simplifying we obtain,
\begin{equation}
\tilde{C}_{1} = 
\dfrac{\kappa_{a}\gamma\, d}{\kappa_{a} - \kappa_{e}}
\left(
\dfrac{1-\beta^{n+1}}{1 - \beta}
\right)
\end{equation}
%---
\end{itemize}
\end{proof}

\subsection{Proof of Theorem \ref{Th:multidose_solution_general}.} \label{AP: generalizedarbitary}

\begin{proof}
    
We prove this theorem via two separate propositions.

%%%%%%%%%%%%%%%%
\begin{proposition}[Recursive Formula for the Generalized Bateman Function]\label{Teo:Formulas_recursivas}
Provided that $\kappa_{a} \neq \kappa_{e}$ then given any arbitrary dosage schedule given by the sequence $\left\{(d_{n}, \tau_{n})\right\}_{n=1}^{\infty}$, the solution to the system of dynamic equations presented in \eqref{Eq:General IVP} is:

%--------
\begin{equation}
\begin{split}
y(t) &= \sum \limits_{n=1}^\infty \mathbbm{1}[t \in I_{n}] \overset{(n)}{y}(t)
\\
x(t) &= \sum \limits_{n=1}^\infty \mathbbm{1}[t \in I_{n}] \overset{(n)}{x}(t)
\end{split}
\end{equation}

where

\begin{equation}
\begin{split}
\stackrel{(n)}{y}(t) &=
\left(
\stackrel{(n-1)}{\mathbf{Rem}_{y}} + d_{n}
\right)
e^{-\kappa_{a}(t - t_{n-1})}, \qquad t \in I_{n} = [t_{n-1}, t_{n}]
\\
\stackrel{(n)}{x}(t) &=
C_{1}(n) \, e^{-\kappa_{e}(t - t_{n-1})}
-
C_{2}(n) \, e^{-\kappa_{a}(t - t_{n-1})}, 
\qquad t \in I_{n} = [t_{n-1}, t_{n}]
\\
C_{1}(n) &= \dfrac{\kappa_{a}\, \gamma }{\kappa_{a} - \kappa_{e}}
\left(
\stackrel{(n-1)}{\mathbf{Rem}_{y}} + d_{n}
\right)
+
\stackrel{(n-1)}{\mathbf{Rem}_{x}}
\\
C_{2}(n) &= \dfrac{\kappa_{a}\, \gamma\ }{\kappa_{a} - \kappa_{e}}
\left(
\stackrel{(n-1)}{\mathbf{Rem}_{y}} + d_{n}
\right)
\\
\text{where} \,\, \stackrel{(n)}{\mathbf{Rem}_{y}} &=
\stackrel{(n)}{y}(t_{n})
\quad 
\text{is the remainder for} \,\, y(t)\,\,
\text{at the end of interval}\,\, I_{n},
\\
\text{and} \,\, \stackrel{(n)}{\mathbf{Rem}_{x}} &=
\stackrel{(n)}{x}(t_{n})
\quad 
\text{is the remainder for} \,\, x(t)\,\,
\text{at the end of interval}\,\, I_{n}
\\
\end{split}    
\end{equation}
\end{proposition}
%%%%%%%%%%%%%%%%%%
%---
\begin{proof}
To prove these recurrence formulas, we will solve the initial value problem (IVP) \eqref{Eq:General IVP} in a generic $n$-th interval, $I_{n} = [t_{n-1}, t_{n}]$, taking into account the initial and multiplicity conditions.\\ 

Denote as $\stackrel{(n)}{\mathbf{Rem}_{x}}$ the value of the drug concentration in the blood plasma at the end of the interval, $I_{n}$, that is,
\begin{equation}
\stackrel{(n)}{\mathbf{Rem}_{x}} = \stackrel{(n)}{x}(t_{n})    
\end{equation}
where $\stackrel{(n)}{x}(t)$ is the solution of the IVP with $t \in I_{n}$. Because this quantity is the blood concentration in the body just before the next drug administration, we call it the remainder for $x(t)$ in the $n-th$ interval. Likewise, we denote as $\stackrel{(n)}{\mathbf{Rem}_{y}}$ is the amount of drug left in the intestinal tract of the interval at the end of $I_{n}$, that is,
\begin{equation}
\stackrel{(n)}{\mathbf{Rem}_{y}} = \lim\limits_{t\to t_{n}^{-}}\stackrel{(n)}{y}(t)    
\end{equation}
where $\stackrel{(n)}{y}(t)$ is the solution of IVP with $t \in I_{n}$. We call it the remainder for $y(t)$ in the $n-th$ interval.
%-------
\begin{itemize}
%---
\item
The first equation of the system \eqref{Eq:General IVP} is independent of $x(t)$. Therefore its solution is:
\begin{equation}
\begin{split}
\dfrac{d\stackrel{(n)}{y}}{dt} &= -\kappa_{a}\stackrel{(n)}{y} \implies 
\stackrel{(n)}{y}(t) = Ce^{-\kappa_{a} t}
\\
\text{where}\quad C &=  \stackrel{(n)}{y}(t_{n-1}^{+})e^{\kappa_{a}t_{n-1}} \qquad \text{(initial condition)}
\\
\implies
\stackrel{(n)}{y}(t) &=
\stackrel{(n)}{y}[t_{n-1}^{+}]e^{-\kappa_{a}(t - t_{n-1}};\qquad
t\in I_{n}
\\
\therefore 
\quad 
\stackrel{(n)}{y}(t) &=
\left(
\stackrel{(n-1)}{\mathbf{Rem}_{y}} + d_{n}
\right)
e^{-\kappa_{a}(t - t_{n-1})};\qquad
t\in I_{n}
\quad
\text{(multiplicity condition)}
\end{split}    
\end{equation}
%---
\item 
Using the previous solution, it follows that
\begin{equation}
\begin{split}
\dfrac{d\stackrel{(n)}{x}}{dt} &= \kappa_{a}\cdot \gamma  \stackrel{(n)}{y}(t)
-
\kappa_{e}\stackrel{(n)}{x} \implies 
\dfrac{d\stackrel{(n)}{x}}{dt}  + 
\kappa_{e}\stackrel{(n)}{x} =
\kappa_{a}\cdot \gamma  \stackrel{(n)}{y}(t) 
\\ 
\dfrac{d}{dt}\left[e^{\kappa_{e}t}\right] &=
e^{\kappa_{e}t} \kappa_{a}\cdot \gamma 
\stackrel{(n)}{y}(t)
\implies
\int d\left[e^{\kappa_{e}} \stackrel{(n)}{x}(t)\right]
=
\int e^{\kappa_{e}}
\kappa_{a}\cdot \gamma  \left(\stackrel{(n - 1)}{\mathbf{Rem}_y}\right)
e^{-\kappa_{a}(t - t_{n-1})} dt
\\
\implies
\stackrel{(n)}{x}(t) &= \kappa_{a}\cdot\gamma \left(\stackrel{(n-1)}{\mathbf{Rem}_y}\right) e^{\kappa_{a} t_{n-1}}
\left[
\dfrac{e^{-\kappa_{a} t}}{\kappa_{a} - \kappa_{e}} + C e^{-\kappa_{e} t}
\right]
\end{split}    
\end{equation}
Finally, using the initial conditions to find the value of the integration constant $C$ and the multiplicity conditions, we can establish the desired result,
\begin{equation}
\stackrel{(n)}{x}(t) =
\left(\dfrac{\kappa_{a}\cdot \gamma \cdot(\stackrel{(n-1)}{\mathbf{Rem}_{y}} + d_{n})}{\kappa_{a} - \kappa_{e}} + \stackrel{(n - 1)}{\mathbf{Rem}_{x}}\right)e^{-\kappa_{e}(t - t_{n-1})}
-
\left(\dfrac{\kappa_{a}\cdot \gamma \cdot(\stackrel{(n-1)}{\mathbf{Rem}_{y}} + d_{n})}{\kappa_{a} - \kappa_{e}} \right)e^{-\kappa_{a}(t - t_{n-1})}
\end{equation}
\end{itemize}
\end{proof}

\vspace{2mm}

Next, we find explicit formulas for the sequences of the remainders\\

\begin{proposition}[Formulas for the Remainders]
We refer to the amount of drug that remains in the intestinal tract, $\stackrel{(n)}{y}(t_{n})$, and the concentration of the drug that still exists in the bloodstream, $\stackrel{(n)}{x}(t_{n})$, at the end of the $n$-th interval or period between each dose, $I_{n} = [t_{n-1}, t_{n}]$. We denote these as $\displaystyle \stackrel{(n)}{\textbf{Rem}_{y}}$ and $\displaystyle \stackrel{(n)}{\textbf{Rem}_{x}}$ respectively. The following expressions represent these quantities:

\begin{equation}
\begin{split}
\stackrel{(n)}{\textbf{Rem}_{y}} &= \stackrel{(n)}{y}(t_{n}) =
\sum\limits_{i=1}^{n}\prod\limits_{j=i}^{n} d_{i}\alpha_{j};
\qquad \alpha_{s} = e^{-\kappa_{a}t_{s}},
\qquad
s=1,2,3,\dots, n\geq1
\\
\stackrel{(n)}{\textbf{Rem}_{x}} &= \stackrel{(n)}{x}(t_{n}) =
\dfrac{\kappa_{a}\cdot \gamma}{\kappa_{a} - \kappa_{e}}
\left[
\sum\limits_{i=1}^{n}\prod\limits_{j=i}^{n} d_{i}\beta_{j}
- 
\sum\limits_{i=1}^{n}\prod\limits_{j=i}^{n} d_{i}\alpha_{j}
\right];
\, \alpha_{s} = e^{-\kappa_{a}\tau_{s}},\,
\beta_{s} = e^{-\kappa_{e}\tau_{s}},
\,
s=1,2,3,\dots, n\geq1\\
\stackrel{(0)}{\textbf{Rem}_{x}}&=\stackrel{(0)}{\textbf{Rem}_{y}}=0
\end{split}    
\end{equation}
\end{proposition}

%%%%%%%%%%%%%%%%%%%%
%---
\begin{proof}
By Proposition \ref{Teo:Formulas_recursivas}, we know that the general solution can be expressed using the recurrence formula:
%--
\begin{equation}
\begin{split}
\stackrel{(n)}{y}(t) &=
\left( 
\stackrel{(n - 1)}{\textbf{Rem}_{y}} + d_{n}
\right)
e^{-\kappa_{a}(t - t_{n-1})}, 
\quad t \in [t_{n-1}, t_{n}], \quad
\tau_{n} = t_{n} - t_{n-1} 
\\
\stackrel{(n)}{x}(t) &=
C_{1} e^{-\kappa_{e}(t - t_{n-1})} -
C_{2} e^{-\kappa_{a}(t - t_{n-1})}
\quad t \in [t_{n-1}, t_{n}], \quad
\tau_{n} = t_{n} - t_{n-1} 
\\
C_{1} &=
\dfrac{\kappa_{a}\cdot\gamma}{\kappa_{a} - \kappa_{e}}
\left(
\stackrel{(n-1)}{\textbf{Rem}_{y}} +\, d_{n}
\right)
+ \stackrel{(n-1)}{\textbf{Rem}_{x}} 
\\
C_{2} &=
\dfrac{\kappa_{a}\cdot\gamma}{\kappa_{a} - \kappa_{e}}
\left(
\stackrel{(n-1)}{\textbf{Rem}_{y}} +\, d_{n}
\right)
\end{split}
\end{equation}
\begin{itemize}
%%%%%
\item 
The proof of these formulas can be found using mathematical induction. Although we will not demonstrate the assertion in this manner, we will outline the steps of the recurrence until the behavior pattern becomes evident.
%---
\begin{equation}
\begin{split}
\stackrel{(0)}{\textbf{Rem}_{y}} &= \stackrel{(0)}{y} (0) = 0
\\
\stackrel{(1)}{\textbf{Rem}_{y}} &=
\stackrel{(1)}{y} (t_{1}) = d_{1} \, \alpha_{1}
\\
\stackrel{(2)}{\textbf{Rem}_{y}} &=
\left(\stackrel{(1)}{\textbf{Rem}_{y}} + d_{2}\right)e^{-\kappa_{a}(t_{2} - t_{1})}
\\
\stackrel{(2)}{\textbf{Rem}_{y}} &=
\left(d_{1}\alpha_{1} + d_{2}\right)\alpha_{2}
\\
\stackrel{(2)}{\textbf{Rem}_{y}} &=
d_{1}\alpha_{1}\alpha_{2} + d_{2}\alpha_{2}
\\
\stackrel{(3)}{\textbf{Rem}_{y}} &=
(d_{1}\alpha_{1}\alpha_{2} + d_{2}\alpha_{2} + d_{3})\alpha_{3}
\\
\stackrel{(3)}{\textbf{Rem}_{y}} &=
d_{1}\alpha_{1}\alpha_{2}\alpha_{3} + d_{2}\alpha_{2}\alpha_{3} + d_{3}\alpha_{3}
\\
\stackrel{(4)}{\textbf{Rem}_{y}} &=
d_{1}\alpha_{1}\alpha_{2}\alpha_{3}\alpha_{4} + d_{2}\alpha_{2}\alpha_{3}\alpha_{4} + d_{3}\alpha_{3}\alpha_{4} + d_{4}\alpha_{4}
\\
\vdots & \qquad \vdots
\\
\stackrel{(n)}{\textbf{Rem}_{y}} &=
\sum\limits_{i=1}^{n}\prod\limits_{j=i}^{n}
d_{i}\alpha_{j}
\end{split}    
\end{equation}

%%%%%
\item 
For the remainder of $\stackrel{(n)}{x}(t)$:
\begin{equation}
\begin{split}
\stackrel{(0)}{\textbf{Rem}_{x}} &=
\stackrel{(0)}{x} (0) = 0
\\
\stackrel{(1)}{\textbf{Rem}_{x}} &=
\stackrel{(1)}{x}(t_{1}) =
C_{1}e^{-\kappa_{e}\tau_{1}} -
C_{2}e^{-\kappa_{a}\tau_{1}}
\\
\stackrel{(1)}{\textbf{Rem}_{x}} &=
\left[\dfrac{\kappa_{a}\cdot \gamma}{\kappa_{a} - \kappa_{e}}\left(\cancelto{0}{\stackrel{(0)}{\textbf{Rem}_{y}}} + d_{1}\right) + \cancelto{0}{\stackrel{(0)}{\textbf{Rem}_{x}}}\right]e^{-\kappa_{e}\tau_{1}} 
-
\left[\dfrac{\kappa_{a}\cdot \gamma}{\kappa_{a} - \kappa_{e}}\left(\cancelto{0}{\stackrel{(0)}{\textbf{Rem}_{y}}} + d_{1}\right) \right]e^{-\kappa_{a}\tau_{1}}
\\
\stackrel{(1)}{\textbf{Rem}_{x}} &=
\dfrac{\kappa_{a}\cdot \gamma}{\kappa_{a} - \kappa_{e}}
\left(d_{1}\beta_{1} - d_{1}\alpha_{1}\right)
\\
\stackrel{(2)}{\textbf{Rem}_{x}} &=
\stackrel{(2)}{x}(t_{2}) =
C_{1}e^{-\kappa_{e}\tau_{2}} -
C_{2}e^{-\kappa_{a}\tau_{2}}
\\
\stackrel{(2)}{\textbf{Rem}_{x}} &=
\left[\dfrac{\kappa_{a}\cdot \gamma}{\kappa_{a} - \kappa_{e}}\left(\stackrel{(1)}{\textbf{Rem}_{y}} + d_{2}\right) + \stackrel{(1)}{\textbf{Rem}_{x}}\right]e^{-\kappa_{e}\tau_{2}} 
-
\left[\dfrac{\kappa_{a}\cdot \gamma}{\kappa_{a} - \kappa_{e}}\left(\stackrel{(1)}{\textbf{Rem}_{y}} + d_{2}\right) \right]e^{-\kappa_{a}\tau_{2}}
\\
\stackrel{(2)}{\textbf{Rem}_{x}} &=
\dfrac{\kappa_{a}\cdot \gamma}{\kappa_{a} - \kappa_{e}}
\left[
\left(
d_{1}\beta_{1}\beta_{2} + d_{2}\beta_{2}
\right)
-
\left(
d_{1}\alpha_{1}\alpha_{2} +
d_{2}\alpha_{2}
\right)
\right]
\\
\stackrel{(3)}{\textbf{Rem}_{x}} &=
\left[\dfrac{\kappa_{a}\cdot \gamma}{\kappa_{a} - \kappa_{e}}\left(\stackrel{(2)}{\textbf{Rem}_{y}} + d_{3}\right) + \stackrel{(2)}{\textbf{Rem}_{x}}\right]e^{-\kappa_{e}\tau_{3}} 
-
\left[\dfrac{\kappa_{a}\cdot \gamma}{\kappa_{a} - \kappa_{e}}\left(\stackrel{(2)}{\textbf{Rem}_{y}} + d_{3}\right) \right]e^{-\kappa_{a}\tau_{3}}
\\
\stackrel{(3)}{\textbf{Rem}_{x}} &=
\dfrac{\kappa_{a}\cdot \gamma}{\kappa_{a} - \kappa_{e}}
\left[
\left(
d_{1}\beta_{1}\beta_{2}\beta_{3} + d_{2}\beta_{2}\beta_{3} + d_{3}\beta_{3}
\right)
-
\left(
d_{1}\alpha_{1}\alpha_{2}\alpha_{3} +
d_{2}\alpha_{2}\alpha_{3} + d_{3}\alpha_{3}
\right)
\right]
\\
& \vdots \qquad \qquad \vdots
\\
\stackrel{(n)}{\textbf{Rem}_{x}} &=
\dfrac{\kappa_{a}\cdot \gamma}{\kappa_{a} - \kappa_{e}}
\left[
\sum\limits_{i=1}^{n}
\prod\limits_{j=i}^{n}
d_{i}\beta_{j}
-
\sum\limits_{i=1}^{n}
\prod\limits_{j=i}^{n}
d_{i}\alpha_{j}
\right]
\end{split}
\end{equation}
\end{itemize}    
\end{proof}

\end{proof}

%%%%%%%%%%%%
\subsection{Proof of Proposition \ref{Prop: AUC_In}} \label{AP: Prop3}
\begin{proof}
Since the solution function $\stackrel{(n)}{x}(t)$ in the time interval $I_{n} = [(n-1)\tau, n\tau]$ satisfies the usual Fundamental Theorem of Calculus, we can find a primitive or antiderivative $\stackrel{(n)}{x}(t)$ of the solution to the initial value problem in the theorem, $\stackrel{(n)}{x}(t)$, and then apply the Fundamental Theorem of Calculus.\\ 
\end{proof}

%%%%%%%%%%%%
\subsection{Proof of Proposition \ref{Prop: Bateman_AUC}} \label{AP: Bateman_AUC}
\begin{proof}
Observe that

\begin{equation}
\begin{split}
AUC_{[0, \infty]} 
&=
\lim\limits_{t\to \infty}
\int_{0}^{t}
\stackrel{(1)}{x}(u) \, du
\\
&= 
\lim\limits_{t\to \infty}
\int_{0}^{t}
\left(
C_{1}e^{-\kappa_{e} u}
-
C_{2}e^{-\kappa_{a} u}
\right)
\, du; \qquad 
C_{1} = C_{2} = 
\dfrac{\kappa_{a}\,\gamma\, d}{\kappa_{a} - \kappa_{e}}
\\
&= 
\, C_{1}\, \lim\limits_{t\to \infty}
\left(
\dfrac{1}{\kappa_{e}} 
- e^{-\kappa_{e}\, t}
-
\dfrac{1}{\kappa_{a}} 
+ e^{-\kappa_{a}\, t}
\right)
\\
&=
\dfrac{\kappa_{a}\,\gamma\, d}{\kappa_{a} - \kappa_{e}}
\left(
\dfrac{1}{\kappa_{e}}
-
\dfrac{1}{\kappa_{a}}
\right)
\end{split}    
\end{equation}    
\end{proof}

%%%%%%%%%%%%
\subsection{Proof of Proposition \ref
{Prop: Maximum}} \label{AP: Maximum}
\begin{proof}
Since the solution given by Theorem \eqref{Th:multidose_solution} is smooth in the usual calculus sense, we only need to verify that the Hilger derivative at the endpoint $(n-1)\tau$ is positive.\\

Without loss of generality, let $\kappa_{a} > \kappa_{e}$ as the proof for $\kappa_{e} > \kappa_{a}$ is similar due to the symmetry of the solution concerning these two parameters. 

\begin{equation}
\begin{split}
\kappa_{a} &> \kappa_{e} \implies
\beta > \alpha \implies \kappa_{a} \beta > \kappa_{e}\alpha 
\implies
\kappa_{a} - \kappa_{e} > \kappa_{a}\beta - \kappa_{e}\alpha
\\
\kappa_{a}(1 - \beta) &> \kappa_{e}(1 - \alpha)
\implies
\dfrac{\kappa_{a}(1 - \beta)}{\kappa_{e}(1 - \alpha)} >1
\end{split}    
\end{equation}

With the previous result, we have proven that $C_{2} > C_{1}$. Therefore, $\stackrel{(n)}{t}_{\text{max}} >0$ and $\stackrel{(n)\triangle}{x}((n-1)\tau) > 0$.\\

Lastly, we find the critical points for each period and verify using the second derivative criterion that at $\stackrel{(n)}{t}_{\text{max}}$, the plasma concentration $\stackrel{(n)}{x}(t)$ reaches its maximum value.
\end{proof}

\subsection{Proof of Theorem \ref{T: SSNew}} \label{AP: SS}

\begin{proof}

Notice that for a fixed $t \in I_{n}$ and $n$:

\begin{align*}
\Big |\stackrel{(n)}{x}(t) -\stackrel{(n-1)}{x}(t-\tau) \Big |&=\Big | \beta^{n} \left(
\dfrac{\kappa_{a}\, \gamma \,d}{\kappa_{a} - \kappa_{e}}
\right)e^{-\kappa_{e}(t -n\tau)}
-  \alpha^{n}\left(
\dfrac{\kappa_{a}\, \gamma \, d}{\kappa_{a} - \kappa_{e}}
\right)
e^{-\kappa_{a}(t - n\tau)}\Big |\\
& =  \left(\dfrac{\kappa_{a}\, \gamma \,d}{V |\kappa_{a} - \kappa_{e}|}
\right)\Big |  \beta^{n} e^{-\kappa_{e}(t -n\tau)} -\alpha^{n} e^{-\kappa_{a}(t - n\tau)} \Big | \\
& \leq   \left(\dfrac{\kappa_{a}\, \gamma \,d}{V |\kappa_{a} - \kappa_{e}|}
\right) \left(\beta^{n} +\alpha^{n} \right) 
\end{align*}
    
Hence,

$$\sup \limits_{t \in I_n} \Big |\stackrel{(n)}{x}(t) - \stackrel{(n-1)}{x}(t-\tau) \Big | \leq   \left(\dfrac{\kappa_{a}\, \gamma \,d}{V |\kappa_{a} - \kappa_{e}|}
\right) \left(\beta^{n} +\alpha^{n} \right) $$

So that

$$0 \leq  \lim \limits_{n \to \infty} \left[ \sup \limits_{t \in I_n} \Big |\stackrel{(n)}{x}(t) - \stackrel{(n-1)}{x}(t-\tau) \Big | \right] \leq \lim \limits_{n \to \infty} \left(\dfrac{\kappa_{a}\, \gamma \,d}{V |\kappa_{a} - \kappa_{e}|}
\right) \left(\beta^{n} +\alpha^{n} \right) =0 $$

Thus, by the squeeze theorem:

$$ \lim \limits_{n \to \infty} \left[\sup \limits_{t \in I_n}   \Big |\stackrel{(n)}{x}(t) - \stackrel{(n-1)}{x}(t-\tau) \Big |  \right]=0$$

\end{proof}

%%%%%%%%%%%%
\subsection{Proof of Proposition \ref{Prop: Bounds}} \label{AP: Bounds}
\begin{proof}
We will use the result from the previous proposition and the definition of steady state. In the steady state, the constants satisfy:
\begin{equation}
C_{1} = \frac{\kappa_{a}d\gamma}{\kappa_{a} - \kappa_{e}} 
\frac{1}{1 - \beta};
\quad
C_{2} = \frac{\kappa_{a}d\gamma}{\kappa_{a} - \kappa_{e}} 
\frac{1}{1 - \alpha}
\end{equation}
Therefore, $\frac{C_{2}}{C_{1}} = \frac{1 - \beta}{1 - \alpha}$. By substituting these values into the result of the previous proposition, we obtain the desired result.
Finally, to demonstrate that this quantity is positive, we can use the chain of implications used earlier: $\kappa_{a} > \kappa_{e}$ implies that $\alpha < \beta$, which further implies that $C_{2} > C_{1}$.\\
\end{proof}

%%%%%%%%
\subsection{Proof of Theorem \ref{Th:steady-state_bounds}.} \label{AP: generalized2}
%---
\begin{proof}

To prove the main statement, we begin proving the following preliminary proposition

\begin{proposition}\label{Th:notation_simplification}
Consider the following short-hand notations:
\begin{equation}
\begin{split}
p_{1} &= -\dfrac{\kappa_{a}}{\kappa_{a} - \kappa_{e}},
\qquad
p_{2} = -\dfrac{\kappa_{e}}{\kappa_{a} - \kappa_{e}},
\\
p_{3} &= 
\left(
\dfrac{\kappa_{a}}{\kappa_{e}}
\right)^{p_{1}},
\qquad
p_{4} = 
\left(
\dfrac{\kappa_{a}}{\kappa_{e}}
\right)^{p_{2}},
\\
z &= 1 - \alpha,
\qquad
w = 1 - \beta.
\end{split}    
\end{equation}
If $\kappa_{a} > \kappa_{e} > 0$, then it follows that
%---
\begin{equation}
\begin{split}
& p_{2} > p_{1}, \qquad p_{4} > p_{3},
\qquad p_{2} - p_{1} = 1 \implies
\begin{cases}
p_{1} = p_{2} - 1
\\
p_{2} = p_{1} + 1
\end{cases}
\\
&\beta > \alpha > 0
\implies
z > w > 0 \implies
\begin{cases}
z' = \frac{d z}{d\, \tau} = \kappa_{a} \alpha
\\[4pt]
w' = \frac{d w}{d\, \tau} = \kappa_{e} \beta
\end{cases}
\\
&\text{and}\quad
\dfrac{w^{p_{1}}}{z^{P_{2}}} > 0
\end{split}   
\end{equation}
\end{proposition}
\begin{proof}
Direct calculations based on the hypothesis and definition of new variables.\\  

\end{proof}

By substituting the new notation defined in Proposition \ref{Th:notation_simplification} and simplifying, the formulas for the minimum and maximum plasma concentrations in the steady-state, $\underline{SS}$ and $\overline{SS}$, are as follows:
%%%%
\begin{equation}\nonumber
\begin{split}
&
\underline{SS}(d, \tau) =
\dfrac{\kappa_{a}\, d\, \gamma}{\kappa_{a} - \kappa_{e}}
\left[ 
\left(
\dfrac{\beta}{1 - \beta}
\right)
-
\left(
\dfrac{\alpha}{1 - \alpha}
\right)
\right]
\\
&
\underline{SS}(d, \tau) =
\dfrac{\kappa_{a}\, d\, \gamma}{\kappa_{a} - \kappa_{e}}
\left[ 
\dfrac{\beta - \alpha}{(1 - \alpha)(1 - \beta)}
\right]
\\
& 
\overline{SS}(d, \tau) = 
\dfrac{\kappa_{a}\, d\, \gamma}{\kappa_{a} - \kappa_{e}}
\left[ 
\dfrac{1}{1 - \beta}
\left(
\dfrac{\kappa_{a}(1 - \beta)}{\kappa_{e} (1 - \alpha)}
\right)
^{
-\dfrac{\kappa_{e}}{\kappa_{a} - \kappa_{e}}
}
-
\dfrac{1}{1 - \alpha}
\left(
\dfrac{\kappa_{a}(1 - \beta)}{\kappa_{e} (1 - \alpha)}
\right)
^{
-\dfrac{\kappa_{a}}{\kappa_{a} - \kappa_{e}}
}
\right] 
\\
& \overline{SS}(d, \tau) =
\dfrac{\kappa_{a}\, d\, \gamma}{\kappa_{a} - \kappa_{e}}
\left(
\dfrac{w^{p_{1}}}{z^{p_{2}}}
\right)
\left(
p_{4} - p_{3}
\right)
\end{split}    
\end{equation}

%--------
\begin{enumerate}
\item [(a)]
Note that both functions are positive because the factors composing them are positive. Additionally, observe that both functions share the same first factor, which is directly proportional to $d$. Therefore, both partial derivatives of the mentioned functions with respect to $d$ are positive, indicating that both functions are increasing with respect to $d$.\\
%%%%%
%%%%%%%
\item [(b)]
By taking partial derivatives with respect to $\tau$ of the expressions obtained in the previous item, we have:

\begin{equation}
\begin{split}
\dfrac{\partial}{\partial \tau}\underline{SS} (d, \tau) &= A
\left[
\dfrac{\kappa_{a}\alpha}{(1 - \alpha)^{2}}
-
\dfrac{\kappa_{e}\beta}{(1 - \beta)^{2}}
\right]
\\
\dfrac{\partial}{\partial \tau} \overline{SS} (d, \tau) &=\dfrac{\partial \overline{SS}}{\partial \tau} =
A (p_{4} - p_{3})
\left(
\dfrac{w^{p_{1}}}{z^{p_{2}}}
\right)
\left(
\dfrac{\kappa_{a}\kappa_{e}}{\kappa_{a} - \kappa_{e}}
\right)
\left[ 
\dfrac{\alpha - \beta}{(1 - \alpha)(1 - \beta)}
\right]
\end{split}    
\end{equation}

Using the results from Proposition \ref{Th:notation_simplification}, it is evident that the function $\overline{SS}$ is decreasing with respect to $\tau$.\\

To demonstrate that the function $\underline{SS}$ is also decreasing with respect to $\tau$, we define a new function $f(x)$ in such a way that the second factor can be rewritten as $f(a) - f(b)$,

\begin{equation}
\begin{split}
f(x) &:= \dfrac{x e^{-x\, \tau}}{(1 - e^{-x\, \tau})^2}
\quad  \text{for} \quad x > 0, \quad \tau > 0.
\\
f(\kappa_{a}) &= \dfrac{\kappa_{a} \alpha}{(1 - \alpha)^2},
\qquad 
f(\kappa_{e}) = \dfrac{\kappa_{e} \beta}{(1 - \beta)^2}
\end{split}
\end{equation}

The function $f(x)$ defined above has a derivative given by,

\begin{equation}
f\,'(x) = 
\dfrac{e^{\tau \, x}
\left[ 
(1 - \tau \, x)e^{\tau \, x 
}
-
(1 + \tau \, x)
\right]}{(e^{\tau \, x} - 1)^3} =
\dfrac{e^{\tau\, x}\left[h_{1}(x) - h_{2}(x) \right]}{(e^{\tau \, x} - 1)^{3}}
\end{equation}

Note that $f\,'(x)$ is negative because $h_{1}(x) = (1 + \tau \, x)e^{\tau \, x} < h_{2}(x) = (1 + \tau \, x)$ due to $\lim\limits_{x\to 0} h_{1}(x) = \lim\limits_{x\to 0} h_{2}(x) = 1$ and additionally, $h_{2}\, ' (x) = \tau >0$ and $h_{1}\, '(x) = -\tau^{2} \, x e^{\tau \, x} < 0$ for all $x>0$.\\

This implies that $f(x)$ is decreasing in the interval $x > 0$. Thus, for $\kappa_{a} > \kappa_{e} > 0$ and $\tau > 0$, we can conclude that $p_{5} = (\kappa_{a} \, \alpha)/(1 - \alpha)^{2} < p_{6} = (\kappa_{e} \, \beta)/(1 - \beta)^{2}$.
Thus, $\underline{SS}$ is decreasing with respect to $\tau$.
\end{enumerate}
\end{proof}

\subsection{Proof of Theorem \ref{Th:Therapeutic_range_width}} \label{AP: TWidth}
%---
\begin{proof}
Let $\ell = \overline{SS} - \underline{SS}$. By using the definitions and results from Proposition \ref{Th:notation_simplification}, we have:

\begin{equation}
\begin{split}
\ell = \overline{SS} - \underline{SS} = 
\dfrac{\kappa_{a}\, d\, \gamma}{\kappa_{a} - \kappa_{e}}
\left[ 
(p_{4} - p_{3})
\left(
\dfrac{w^{p_{1}}}{z^{p_{2}}}
\right)
+
\left(
\dfrac{z - w}{z\, w}
\right)
\right]
\end{split}    
\end{equation}

Note that $\ell > 0$ for all $\kappa_{a} > \kappa_{e} > 0$, $d>0$, and $\tau >0$.\\
To prove the second part, we again use the definitions from Proposition \ref{Th:notation_simplification}. Note that if $\tau \to \infty$, then $\alpha \to 0$, $\beta \to 0$, $z \to 1$, and $w \to 1$. We obtain the desired result by substituting the original variables $p_{i}$, $w$, and $z$.\\
\end{proof}

%%%%%%%%%%%%
\subsection{Proof of Proposition \ref{Th: equality_AUCs}} \label{AP: equality_AUCs}
\begin{proof}
From the definition of steady state, it follows that
\begin{equation}
\text{AUC}^{\,s.s.} = \lim\limits_{n\to \infty} \text{AUC}_{I_{n}}    
\end{equation}
Under the assumption that $\kappa_{a} > \kappa_{e} >0$, $\alpha = e^{-\kappa_{a}\tau}$ and $\beta = e^{-\kappa_{e}\tau}$, then $0 < \alpha < \beta < 1$ from which it can be contended that,
\begin{equation}
\lim\limits_{n\to \infty} \alpha^{n} 
=
\lim\limits_{n\to \infty} \beta^{n} = 0
\end{equation}
Now, following Proposition \ref{Prop: AUC_In} and Proposition \ref{Prop: Bateman_AUC} we can conclude that,
\begin{equation}
\text{AUC}^{\,s.s.} =
\lim\limits_{n\to\infty}
\dfrac{\kappa_{a}\,\gamma\, d}{\kappa_{a} - \kappa_{e}}
\left[
\left(
\dfrac{1 - \beta^{n}}{\kappa_{e}}
\right)
-
\left(
\dfrac{1 - \alpha^{n}}{\kappa_{a}}
\right)
\right]
= \text{AUC}_{[0, \infty)}
\end{equation}
\end{proof}

\subsection{Proof of Theorem \ref{T: Existence}} \label{AP: Existence}

\begin{proof}

For any $\bar{R}>\underline{R}>0$, we can choose target levels $(\underline{SS}^*,\overline{SS}^*)$ such that $\underline{R}<\underline{SS}^*<\overline{SS}^*<\overline{R}$. We know establish the solutions to the 
system of nonlinear equations given by,

\begin{equation}
\begin{split}
\underline{SS} = \underline{SS}(\kappa_{a}, \kappa_{e}, \gamma, V, d, \tau) 
\\
\overline{SS} = \overline{SS}(\kappa_{a}, \kappa_{e}, \gamma, V, d, \tau) 
\end{split}    
\end{equation}

Furthermore, for a fixed vector of physiological parameters $(\kappa_a,\kappa_e,\gamma)$, we can further view this problem as 

\begin{equation}
\begin{split}
\underline{SS} = \underline{SS}(d, \tau) 
\\
\overline{SS} = \overline{SS}(d, \tau) 
\end{split}    
\end{equation}

Within this logical framework, we must choose a pair  $(d^*,\tau^*)$ such that

$$\underline{SS}(d^*,\tau^*)=\dfrac{\kappa_{a}\,d\,\gamma}{\kappa_{a} - \kappa_{e}} \left[
\left(
\dfrac{\beta(\tau) }{1 -\beta(\tau)}
\right)
-
\left(
\dfrac{\alpha(\tau)}{1 -\alpha(\tau)}
\right) \right]=\dfrac{\kappa_{a}\,d\,\gamma}{\kappa_{a} - \kappa_{e}} \Psi(\tau)=\underline{SS}^*$$

Hence, all solutions must satisfy

$$d^*(\tau)=\dfrac{\kappa_{a} - \kappa_{e}}{\kappa_{a}\,\gamma} \dfrac{\underline{SS}^*}{\Psi(\tau)} $$

Restricted to the curve, the upper limit can be seen as a function of $\tau$ as follows

\begin{align*}
\overline{\text{SS}}(d^*(\tau),\tau) &= 
\frac{\kappa_{a}d^*(\tau)\gamma}{\kappa_{a} - \kappa_{e}}
\left[
\frac{1}{1 - \beta}
\left(
\frac{\kappa_{a}(1-\beta)}{\kappa_{e}(1-\alpha)}
\right)^{-\frac{\kappa_{e}}{\kappa_{a} - \kappa_{e}}}
-
\frac{1}{1 - \alpha}
\left(
\frac{\kappa_{a}(1-\beta)}{\kappa_{e}(1-\alpha)}
\right)^{-\frac{\kappa_{a}}{\kappa_{a} - \kappa_{e}}}
\right]\\
&=\frac{\kappa_{a}d^*(\tau)\gamma}{\kappa_{a} - \kappa_{e}} \Phi(\tau)\\
&=\dfrac{\Phi(\tau)}{\Psi(\tau)}\underline{SS}^*
\end{align*}

Thus, it suffices to show there is a $\tau>0$ satisfying the equation

$$\overline{\text{SS}}(d^*(\tau),\tau)=\dfrac{\Phi(\tau)}{\Psi(\tau)}\underline{SS}^*=\overline{SS}^*$$

Equivalently, the question boils down to determining

$$\dfrac{\overline{SS}^*}{\underline{SS}^*} \overset{?}{\in} \text{Range}(f(\tau)); \quad \text{ where } f(\tau)=\dfrac{\Phi(\tau)}{\Psi(\tau)}$$

\begin{proposition}[Range of the quotient]\label{Th:Range}
Let \(\kappa_{a} > \kappa_{e}\). Then the function $f(\tau)$ defined as

\begin{equation}
\begin{split}
f(\tau) &= \dfrac{\Phi(\tau;\kappa_{a}, \kappa_{e}, \gamma, V)}{\Psi(\tau;\kappa_{a}, \kappa_{e}, \gamma, V)}
\\
\Psi(\tau) &= \dfrac{\beta}{ 1 - \beta}
- \dfrac{\alpha}{1 - \alpha}
\\
\Phi(\tau) &= \dfrac{1}{ 1 - \beta}
\left(
\dfrac{\kappa_{a}(1 - \beta)}{\kappa_{e}(1 - \alpha)}
\right)^{p_{2}}
- \dfrac{1}{1 - \alpha}
\left(
\dfrac{\kappa_{a}(1 - \beta)}{\kappa_{e}(1 - \alpha)}
\right)^{p_{1}}
\\
\text{where} &\quad
\alpha = e^{-\kappa_{a}\,\tau}
\qquad
\beta = e^{-\kappa_{e}\,\tau}
\qquad
p_{1} = -\dfrac{\kappa_{a}}{\kappa_{a} - \kappa_{e}}
\qquad
p_{2} = -\dfrac{\kappa_{e}}{\kappa_{a} - \kappa_{e}}
\end{split}
\end{equation}
is well-defined in the domain $(0, \infty)$, has a range of $(1, \infty)$ and is a continuous function.

\end{proposition}

%%%%%%%%%%
\begin{proof}
To demonstrate that the range of $f(\tau)$is $(1, \infty)$ in the domain $(0, \infty)$, we will consider the following auxiliary functions:

\begin{equation}
\begin{split}
g(\tau) &=
\dfrac{\kappa_{a}(1 - \beta)}
{\kappa_{e}(1 - \alpha)}
\\
f_{1}(\tau) &=
\dfrac{1}{\Psi(\tau)}
\left(
\dfrac{1}{1 - \alpha}
\right)
g(\tau)^{p_{1}} =
\left(
\dfrac{1 - \beta}{\beta - \alpha}
\right)
f(\tau)^{p_{1}}
\\
f_{2}(\tau) &=
\dfrac{1}{\Psi(\tau)}
\left(
\dfrac{1}{1 - \beta}
\right)
g(\tau)^{p_{2}} =
\left(
\dfrac{1 - \alpha}{\beta - \alpha}
\right)
f(\tau)^{p_{2}}
\\
f(\tau) &= f_{2}(\tau) - f_{1}(\tau) 
\end{split}    
\end{equation}

\begin{enumerate}[a)]
\item
Following L'Hôpital's rule, we have:
\begin{equation}
\lim\limits_{\tau\to 0} g(\tau) = 
\left(
\dfrac{\kappa_{a}}{\kappa_{e}}
\right)
\lim\limits_{\tau\to 0}
\left(
\dfrac{1 - \beta}{1 - \alpha}
\right)
=
\left(
\dfrac{\kappa_{a}}{\kappa_{e}}
\right)
\lim\limits_{\tau\to 0}
\left(
\dfrac{\kappa_{e}\beta}{\kappa_{a} \alpha}
\right) = 1
\end{equation}

Now,
\begin{equation}
\begin{split}
\lim\limits_{\tau\to 0} f_{2}(\tau) = 
\lim\limits_{\tau\to 0}
\left(
\dfrac{1 - \alpha}{\beta - \alpha}
\right)
\left(
\lim\limits_{\tau\to 0} f(\tau)
\right)^{p_{2}}
=
\lim\limits_{\tau\to 0}
\left(
\dfrac{\kappa_{a}\alpha}
{\kappa_{e}\beta - \kappa_{a}\alpha}
\right)
\left(
\lim\limits_{\tau\to 0} f(\tau)
\right)^{p_{2}} = - p_{1}
\end{split}    
\end{equation}

and similarly,
\begin{equation}
\begin{split}
\lim\limits_{\tau\to 0} f_{1}(\tau) = 
\lim\limits_{\tau\to 0}
\left(
\dfrac{1 - \beta}{\beta - \alpha}
\right)
\left(
\lim\limits_{\tau\to 0} f(\tau)
\right)^{p_{1}}
=
\lim\limits_{\tau\to 0}
\left(
\dfrac{\kappa_{e}\beta}
{\kappa_{e}\beta - \kappa_{a}\alpha}
\right)
\left(
\lim\limits_{\tau\to 0} f(\tau)
\right)^{p_{1}} = - p_{2}
\end{split}    
\end{equation}

Therefore,
\begin{equation}
\lim\limits_{\tau\to 0} f(\tau) =
\lim\limits_{\tau\to 0} f_{1}(\tau)
-
\lim\limits_{\tau\to 0} f_{2}(\tau) = 
p_{2} - p_{1} = 1
\end{equation}

\item  %b
Since $f(\tau) = \dfrac{\Phi(\tau)}{\Psi(\tau)}$, we have:
\begin{equation}
\begin{split}
\lim\limits_{\tau\to\infty}
\Psi(\tau) &= 0
\\
\lim\limits_{\tau\to\infty}
\Phi(\tau) &= 
\left(
\dfrac{\kappa_{a}}{\kappa_{e}}
\right)^{p_{2}}
-
\left(
\dfrac{\kappa_{a}}{\kappa_{e}}
\right)^{p_{1}} > 0
\qquad
\text{because}
\quad
\kappa_{a} > \kappa_{e}
\quad
\text{and}
\quad
p_{2} > p_{1}
\end{split}    
\end{equation}

Therefore,
\begin{equation}
\lim\limits_{\tau\to\infty}
f(\tau) = 
\dfrac{\lim\limits_{\tau\to\infty} \Phi(\tau)}{\lim\limits_{\tau\to\infty} \Psi(\tau)} = \infty
\end{equation}

\item %c
Because $f(\tau)$ is the ratio of analytic functions, $\Phi(\tau)$ and $\Psi(\tau)$, defined on $(0, \infty)$, then $f(\tau)$ is also continuous on $(0, \infty)$.\\
\end{enumerate}

Given $\lim \limits_{\tau \to 0} f(\tau)=1$, $\lim \limits_{\tau \to \infty} f(\tau)=\infty$ and $f(\tau)$ is continuous on $(0,\infty)$, it follows by the Intermediate Value theorem that the range of $f(\tau)$ is $(1,\infty)$.\\

\end{proof}
\vspace{2mm}

Following Proposition \ref{Th:Range}, for any $(\underline{SS}^*,\overline{SS}^*)$ such that $\underline{R}<\underline{SS}^*<\overline{SS}^*<\overline{R}$, there exists $\tau^*=\tau^*(\underline{SS}^*,\overline{SS}^*)$ such that

$$f(\tau^*)=\dfrac{\overline{SS}^*}{\underline{SS}^*}$$

Thus the pair $(d^*,\tau^*)=(d^*(\tau^*(\underline{SS}^*,\overline{SS}^*)), \tau^*(\underline{SS}^*,\overline{SS}^*))$ is such that

\begin{equation}
\begin{split}
\underline{SS}^* = \underline{SS}(d^*, \tau^*) 
\\
\overline{SS}^* = \overline{SS}(d^*, \tau^*) 
\end{split}    
\end{equation}
\vspace{2mm}

which means that $(d^*,\tau^*) \in  \mathcal{E}(\underline{R},\bar{R};\kappa_a,\kappa_e,\gamma)$, so that $ \mathcal{E}(\underline{R},\bar{R};\kappa_a,\kappa_e,\gamma) \neq \emptyset$.

\end{proof}

Furthermore, we can prove the effective dose is unique for a given target level $(\underline{SS}^*,\overline{SS}^*)$, since the map $f(\tau)$ introduced above is also injective:\\

\begin{proposition}[The function $f(\tau)$ is strictly increasing]\label{Prop: Increasing}
The function $f(\tau): (0,\infty) \longrightarrow (1, \infty)$ is strictly increasing, i.e., $f(\tau)$ is bijective.
\end{proposition}
%---
\begin{proof}

In Theorem \ref{Th:Range} we proved that $f(\tau)$, with a domain of $[0, \infty)$ and a range of $[1, \infty)$, is continuous and satisfies: $\lim\limits_{\tau \to 0} f(\tau) = 1$ and $\lim\limits_{\tau \to \infty} f(\tau) = \infty$. 
To prove that $f(\tau)$ is bijective, we will show that it is strictly increasing by arguing that its derivative is positive, $f\,'(\tau) >0$, for all $\tau > 0$.\\

Departing from the definition of $f(\tau)$, we can derive the following simplified expression of its derivative:\\
%----
\begin{equation}
\begin{split}
f\,'(\tau) &= 
\left[
\left(
\dfrac{\kappa_{a}}{\kappa_{e}}
\right)^{p_{2}}
-
\left(
\dfrac{\kappa_{a}}{\kappa_{e}}
\right)^{p_{1}}
\right] 
\left[
\dfrac{(1 - \beta)^{p_{1}}}{(1 - \alpha)^{p_{2}}}
\right]
\left[
\dfrac{\kappa_{e} \beta (1 - \alpha)^{2}}{(\beta - \alpha)^{2}}
-
\dfrac{\kappa_{a} \alpha (1 - \beta)^{2}}{(\beta - \alpha)^{2}}
-
\dfrac{\kappa_{a} \, \kappa_{e}}{\kappa_{a} - \kappa_{e}}
\right]
\end{split}    
\end{equation}

Note that the first two of the three factors that compose the derivative are positive. It remains to prove that the third factor is also positive. To analyze the third factor more easily, we will use the following notation: let $a = \kappa_{a}$, $b = \kappa_{e}$, $z = 1 -\alpha$, $w = 1 -\beta$. Since $a > b > 0$ and $\tau >0$, we have $0 < w < z < 1$. With this new notation, the derivative can be expressed as:
\begin{equation}
\begin{split}
f\,'(\tau) &=
\left[
\left(
\dfrac{a}{b}
\right)^{p_{2}}
-
\left(
\dfrac{a}{b}
\right)^{p_{1}}
\right] 
\left[
\dfrac{w^{p_{1}}}{z^{p_{2}}}
\right]
\left[
\dfrac
{
a^{2} w^{2} z - a^{2} w^{2} - a b w^{2} z - a b w z^{2} + 2 a b w z + b^{2} w z^{2} - b^{2} z^{2}
}
{
(a - b)(z - w)^{2}
}
\right]
\\
f\,'(\tau) &=
\left[
\left(
\dfrac{a}{b}
\right)^{p_{2}}
-
\left(
\dfrac{a}{b}
\right)^{p_{1}}
\right] 
\left[
\dfrac{w^{p_{1}}}{z^{p_{2}}}
\right]
\left[
\dfrac
{
bz(1 - w)(aw - bz) -
aw(1 - z)(aw - bz)
}
{
(a - b)(z - w)^{2}
}
\right]
\end{split}
\end{equation}
Let's denote by $\text{num}$ the numerator of the third factor of the derivative of $f(\tau)$:
\begin{equation}
\text{num} = 
bz(1 - w)(aw - bz) -
aw(1 - z)(aw - bz)
\end{equation}
Let's prove that $\text{num} = bz(1 - w)(aw - bz) -
aw(1 - z)(aw - bz) >0$. This can be accomplished in two steps:
\\

\textbf{Step 1:} \underline{Prove that $aw - bz >0$:} Consider the function $\displaystyle g(x) = \dfrac{xe^{x\tau}}{e^{x\tau} -1}$ and prove that $g(\cdot)$ is increasing. The derivative, $\displaystyle g\,'(x) = \dfrac{e^{x\tau}\left[e^{x\tau} - (1 + x\tau)\right]}{(e^{x\tau} - 1)^{2}}$, is positive because the second term in the numerator dominates (is greater than) the first term in the numerator. This implies that $g(a) > g(b)$, which implies that indeed $\displaystyle aw - bz >0$.\\
\\

\textbf{Step 2:} \underline{Prove that $\displaystyle \text{num} = 
bz(1 - w)(aw - bz) -
aw(1 - z)(aw - bz)$:} 

Consider the function $\displaystyle h(x) = \dfrac{x}{e^{x\tau} - 1}$. 

This new function is decreasing because its derivative is negative. This is because $\displaystyle h\,'(x) = \dfrac{e^{x\tau} - (1 + x\tau e^{x\tau})}{(e^{x\tau} - 1)^{2}}$, and the second term in the numerator dominates the first term in the numerator. Therefore, $h(a) < h(b)$, which implies that: $\displaystyle \text{num} = 
bz(1 - w)(aw - bz) -
aw(1 - z)(aw - bz) >0$.\\

Hence, we can conclude that for all $\tau >0$ and for any choice of $\kappa_{a} > \kappa_{e} >0$, the function $f(\tau)$ is strictly increasing on $[0, \infty)$ and thus also bijective.

\end{proof}

An immediate corollary of Proposition \ref{Prop: Increasing} is that the effective doses can be locally expressed  as functions of the biological parameters and the target concentration levels:

\[
\begin{split}
\tau^*&=\tau^*(\underline{\text{SS}}^*,\overline{\text{SS}}^*,\kappa_a,\kappa_e,\gamma)\\
d^*&=d^*(\underline{\text{SS}}^*,\overline{\text{SS}}^*,\kappa_a,\kappa_e,\gamma)
\end{split} 
\]

\newpage
%% ====== SECTION 8 ===== %%
\section{Appendix - Omitted proofs - Other applications}

\subsection{Proof of Theorem \ref{Th:multidose_bolus}}\label{AP: BolusMulti}

\begin{proof}

To establish the validity of the theorem, we prove an equivalent proposition

\begin{proposition}
Given any arbitrary dosage schedule given by the sequence $\left\{(\delta_{n}, \tau_{n})\right\}_{n=1}^{\infty}$, the solution to the system of dynamic equations presented in \eqref{Eq:SEC4_EXAMPLE_1} is:
%% ================= %

\begin{equation}\label{Eq:SEC4_EXAMPLE_1_SOL_REM}
\begin{split}
x(t) &= \sum \limits_{n=1}^{\infty} \mathbbm{1}[t \in I_{n}] \overset{\scriptstyle{(n)}}{x}(t); 
\quad 
\overset{\scriptstyle{(n)}}{x}(t) =  
\left(
\mathbf{\overset{(n-1)}{Rem_{x}}}
 + \delta_{n}
\right)
e^{-\kappa_{e}(t - t_{n-1})},
\quad
\beta_{s} = e^{-\kappa_{e}\tau_{s}} 
\\
\text{where} &
\\
\mathbf{\overset{(n)}{Rem_{x}}} &=
\sum\limits_{i=1}^{n}
\prod\limits_{j=i}^{n}\delta_{i}\beta_{j}
\quad
\text{and}
\quad
\mathbf{\overset{(0)}{Rem_{x}}} = 0
\end{split}
\end{equation}
\end{proposition}

\begin{proof}

We will prove the theorem using mathematical induction.
\begin{enumerate}[(a)]
\item 
\textbf{Base case:}
We verify that the formula holds for $n=1$. Using the initial condition and the continuity condition, we have,
\begin{equation}
\begin{split}
\overset{(1)}{x}(t) &= C_{0}e^{\kappa_{e} t}; 
\qquad t \in [0, t_{1}]
\\
\implies 
C_{0} &= \delta_{1}
\quad
\therefore
\quad
C_{0} = \mathbf{\overset{(0)}{Rem_{x}}} + \delta_{1}
\end{split}    
\end{equation}
\item 
\textbf{Induction hypothesis:}
We assume that the formula holds for $n=k$, that is,
\begin{equation}
\begin{split}
\overset{\scriptstyle{(k)}}{x}(t) &=  
\left(
\mathbf{\overset{(k-1)}{Rem_{x}}}
 + \delta_{k}
\right)
e^{-\kappa_{e}(t - t_{k-1})},
\quad
t \in [t_{k-1}, t_{k}] 
\\
\text{where} &
\\
\mathbf{\overset{(k)}{Rem_{x}}} &=
\sum\limits_{i=1}^{k}
\prod\limits_{j=i}^{k}\delta_{i}\beta_{j}
\end{split}
\end{equation}
\item 
\textbf{Induction step:}
Using the initial condition, the continuity condition, and the induction hypothesis, we have,
\begin{equation}
\begin{split}
\overset{(k+1)}{x}(t) &=
C_{k}e^{\kappa_{e}(t - t_{k})}, 
\qquad
t \in [t_{k}, t_{k+1}]
\\
\implies
C_{k} &= 
\overset{(k)}{x}(t_{k}) + \delta_{k+1}
\\
&= \left(
\mathbf{\overset{(k-1)}{Rem_{x}}} + \delta_{k}
\right)e^{-\kappa_{e}(t_{k} - t_{k-1})}
+ \delta_{k+1}
\\
&= \left(
\mathbf{\overset{(k-1)}{Rem_{x}}} + \delta_{k}
\right)\beta_{k}
+ \delta_{k+1}
\\
&= \left(
\sum\limits_{i=1}^{k-1}
\prod\limits_{j=i}^{k-1}\delta_{i}\beta_{j}
+ \delta_{k}
\right)\beta_{k}
+ \delta_{k+1}
\\
C_{k}
&= 
\sum\limits_{i=1}^{k}
\prod\limits_{j=i}^{k}\delta_{i}\beta_{j}
+ \delta_{k+1}
\\
\therefore
\qquad
\overset{\scriptstyle{(k+1)}}{x}(t) &=  \left(
\mathbf{\overset{(k)}{Rem_{x}}}
 + \delta_{k+1}
\right)
e^{-\kappa_{e}(t - t_{k})},
\quad
t \in [t_{k}, t_{k+1}]
\end{split}
\end{equation}
\end{enumerate}
\end{proof}    
\end{proof}

\subsection{Proof of Theorem \ref{Th:PBFTPK_2}}\label{AP: FATMulti}

\begin{proof}

To establish the validity of the theorem, we prove an equivalent proposition

\begin{proposition}
    
Given any arbitrary dosage schedule given by the sequence $\left\{(\delta_{n}, \tau_{n})\right\}_{n=1}^{\infty}$,  and a sequence of F.A.T. $\{s_n\}_{n=1}^\infty$, the solution to the system of dynamic equations presented in \eqref{Eq:FAT_IVP} is:

\begin{equation}\label{Eq:FAT_IVP_SOL2}
\begin{cases}
\, y(t)= \sum \limits_{n=1}^{\infty} \mathbbm{1}[t \in I_{n}^{1}] \,\,\overset{\scriptstyle{(n)}}{y_{1}}(t), \quad
\overset{\scriptstyle{(n)}}{y_{1}}(t) =  d_{n}
e^{-\kappa_{a}(t - t_{n-1})} 
\\
x(t) = \sum \limits_{n=1}^\infty \mathbbm{1}[t \in I_{n}^{1}] \,\,\overset{\scriptstyle{(n)}}{x_{1}}(t)+\mathbbm{1}[t \in I_{n}^{2}] \,\,\overset{\scriptstyle{(n)}}{x_{2}}(t); \quad 
\begin{cases}
\overset{\scriptstyle{(n)}}{x_{1}}(t) =
C_{1} e^{-\kappa_{e}(t - t_{n-1})}
-
C_{2} e^{-\kappa_{a}(t - t_{n-1})}
\\
\overset{\scriptstyle{(n)}}{x_{2}}(t) =
C_{3} e^{-\kappa_{e}(t - s_{n})}  
\end{cases}
\\
\text{where}
\\[10pt]
I_{n}^{1} = [t_{n-1}, s_{n}],\textbf{(Assimilation Phase)};
\quad 
I_{n}^{2} = [s_{n}, t_{n}],\textbf{(Clearance  Phase)}; \quad n = 1,\dots,\infty
\\[10pt]
C_{1} = 
\dfrac{\kappa_{a}\cdot\gamma}{\kappa_{a} - \kappa_{e}}
d_{n} +
\mathbf{\overset{(n-1)}{Rem^{2}_{x}}}
\\[10pt]
C_{2} = 
\dfrac{\kappa_{a}\cdot\gamma}{\kappa_{a} - \kappa_{e}} d_{n}
\\[15pt]
C_{3} = 
\mathbf{\overset{(n)}{Rem^{1}_{x}}}
\\[10pt]
\quad
\mathbf{\overset{(0)}{Rem^{2}_{x}}} = 0
\end{cases}    
\end{equation}
% \end{theorem}

where,
\begin{equation}\label{Eq:FAT_IVP_REMAINS2}
\begin{split}
\mathbf{\overset{(n)}{Rem^{1}_{x}}} &=
\dfrac{\kappa_{a}\cdot\gamma}{\kappa_{a} - \kappa_{e}}
(B - A)
\left[
\sum\limits_{i = 1}^{n}
\beta^{n-i}\,
d_{i}
\right]
\quad n = 1, 2, \dots,\infty
\\[15pt]
\mathbf{\overset{(n)}{Rem^{2}_{x}}} &=
\dfrac{\beta}{B}
\mathbf{\overset{(n)}{Rem^{1}_{x}}} 
\quad n = 1,\dots,\infty
\\[15pt]
\alpha &= e^{-\kappa_{a} \tau}, 
\quad 0 < \alpha < 1
\quad 
\beta = e^{-\kappa_{e} \tau},
\quad 0 < \beta < 1  
\\
A &= e^{-\kappa_{a} \sigma}, 
\quad 0 < A < 1
\quad 
B = e^{-\kappa_{e} \sigma},
\quad 0 < B < 1  
\\
\alpha^{n} &= e^{-\kappa_{a}t_{n}},
\quad
\beta^{n} = e^{-\kappa_{e}t_{n}},
\quad
\alpha^{n-1} A = e^{-\kappa_{a}s_{n}},
\quad
\beta^{n-1} B = e^{-\kappa_{e}s_{n}}; \quad n = 1,2,\dots,\infty
\\
\alpha^{0} &= \beta^{0}  = 1.    
\end{split}    
\end{equation}
\end{proposition}

\begin{proof}

We first proof  a pair of propositions before dealing with the main result

\vspace{1cm}

%%%%%%%%%%%%%%%%%%%%%%%% PROPOSITION 1
\begin{proposition}\label{Th:Prop_PBFTPK_1}
The Initial Value Problem (IVP),
\begin{equation}
\begin{split}
\begin{cases}
\begin{rcases}
y'(t) = - \kappa_{a} y(t)
\\
x'(t) = \kappa_{a}\gamma y(t) - \kappa_{e} x(t)
\end{rcases}
\quad t \in [t_{n-1}, s_{n}]
\\
\textbf{with initial conditions,}
\\
y(t_{n-1}) = y_{10}
\\
x(t_{n-1}) = x_{10}
\end{cases}    
\end{split}    
\end{equation}
has a unique solution given by,
\begin{equation}
\begin{split}
y(t) &= \left(\dfrac{\kappa_{a}\gamma}{\kappa_{a} - \kappa_{e}} y_{10}\right)
e^{-\kappa_{a}(t - t_{n-1})}
\\[10pt]
x(t) &=
\left(\dfrac{\kappa_{a}\gamma}{\kappa_{a} - \kappa_{e}} y_{10} + x_{10}\right)
e^{-\kappa_{e}(t - t_{n-1})}
-
\left(\dfrac{\kappa_{a}\gamma}{\kappa_{a} - \kappa_{e}} y_{10}\right)
e^{-\kappa_{a}(t - t_{n-1})}
\end{split}    
\end{equation}
\end{proposition}
% ---
\begin{proof}
The matrix $A$ associated with this first-order linear system has eigenvalues and eigenvectors,
\begin{equation}
A =
\begin{bmatrix}
-\kappa_{a} & 0 
\\[10pt]
\kappa_{a}\gamma & -\kappa_{e}
\end{bmatrix}
\longrightarrow
\begin{cases}
\lambda_{1} = -\kappa_{a} 
\rightarrow
\xi_{1} =
\begin{bmatrix}
\dfrac{\kappa_{a} - \kappa_{e}}{\kappa_{a} \gamma}
\\[10pt]
1
\end{bmatrix}
\\[20pt]
\lambda_{2} = -\kappa_{e} 
\rightarrow 
\xi_{2} =
\begin{bmatrix}
0
\\[10pt]
1
\end{bmatrix}
\end{cases}
\end{equation}
Thus, the solution is,
\begin{equation}
\begin{bmatrix}
y(t)
\\[10pt]
x(t)
\end{bmatrix}
=
C_{1}
\begin{bmatrix}
\dfrac{\kappa_{e} - \kappa_{a}}{\kappa_{a} \gamma} 
\\[10pt]
1
\end{bmatrix}
e^{-\kappa_{a}(t - t_{n-1})}
+
C_{2}
\begin{bmatrix}
0 
\\[10pt]
1
\end{bmatrix}
e^{-\kappa_{e}(t - t_{n-1})}
\end{equation}
Replacing the initial conditions at $t = t_{n-1}$, we have,
\begin{equation}
\begin{split}
y(t) &= \left(\dfrac{\kappa_{a}\gamma}{\kappa_{a} - \kappa_{e}} y_{10}\right)
e^{-\kappa_{a}(t - t_{n-1})}
\\[10pt]
x(t) &=
\left(\dfrac{\kappa_{a}\gamma}{\kappa_{a} - \kappa_{e}} y_{10} + x_{10}\right)
e^{-\kappa_{e}(t - t_{n-1})}
-
\left(\dfrac{\kappa_{a}\gamma}{\kappa_{a} - \kappa_{e}} y_{10}\right)
e^{-\kappa_{a}(t - t_{n-1})}
\end{split}    
\end{equation}
\end{proof}

%%%%%%%%%%%%%%%%%%%%%%%% PROPOSITION 2
\begin{proposition}\label{Th:Prop_PBFTPK_2}
The Initial Value Problem (IVP),
\begin{equation}
\begin{split}
\begin{cases}
\begin{rcases}
y'(t) = 0
\\
x'(t) =  - \kappa_{e} x(t)
\end{rcases}
\quad t \in [s_{n}, t_{n}]
\\
\textbf{with initial conditions,}
\\
y(s_{n}) = 0
\\
x(s_{n}) = x_{20}
\end{cases}    
\end{split}    
\end{equation}
has a unique solution given by,
\begin{equation}
\begin{split}
y(t) &= 0
\\[10pt]
x(t) &=
x_{20}
e^{-\kappa_{e}(t - s_{n})}
\end{split}    
\end{equation}
\end{proposition}
% ---
\begin{proof}
The matrix $A$ associated with this first-order linear system has eigenvalues and eigenvectors,
\begin{equation}
A =
\begin{bmatrix}
0 & 0 
\\[10pt]
0 & -\kappa_{e}
\end{bmatrix}
\longrightarrow
\begin{cases}
\lambda_{1} = 0
\rightarrow
\xi_{1} =
\begin{bmatrix}
1
\\[10pt]
0
\end{bmatrix}
\\[20pt]
\lambda_{2} = -\kappa_{e} 
\rightarrow 
\xi_{2} =
\begin{bmatrix}
0
\\[10pt]
1
\end{bmatrix}
\end{cases}
\end{equation}
Thus, the solution is,
\begin{equation}
\begin{bmatrix}
y(t)
\\[10pt]
x(t)
\end{bmatrix}
=
C_{1}
\begin{bmatrix}
1
\\[10pt]
0
\end{bmatrix}
e^{-\kappa_{a}(t - s_{n})}
+
C_{2}
\begin{bmatrix}
0 
\\[10pt]
1
\end{bmatrix}
e^{-\kappa_{e}(t - s_{n})}
\end{equation}
Replacing the initial conditions at $t = s_{n}$, we have,
\begin{equation}
\begin{split}
y(t) &= \left(\dfrac{\kappa_{a}\gamma}{\kappa_{a} - \kappa_{e}} y_{10}\right)
e^{-\kappa_{a}(t - t_{n-1})}
\\[10pt]
x(t) &=
\left(\dfrac{\kappa_{a}\gamma}{\kappa_{a} - \kappa_{e}} y_{10} + x_{10}\right)
e^{-\kappa_{e}(t - t_{n-1})}
-
\left(\dfrac{\kappa_{a}\gamma}{\kappa_{a} - \kappa_{e}} y_{10}\right)
e^{-\kappa_{a}(t - t_{n-1})}
\end{split}    
\end{equation}
\end{proof}

Now we proceed with the main proof:\\

\begin{enumerate}[(a)]
\item 
To prove the formulas for the solutions given in \eqref{Eq:FAT_IVP_SOL2}, we proceed as follows:
\begin{itemize}
\item 
The formula for $\overset{\scriptstyle{(n)}}{y_{1}}(t)$ is obtained by solving the associated equation for $y(t)$ (see Proposition \ref{Th:Prop_PBFTPK_2}) and using the first initial condition and the first multiplicity condition yields the result.\\
\item 
It is clear that from the third multiplicity condition we have $\overset{\scriptstyle{(n)}}{y_{0}}(t) = 0$, which is not explicitly written.\\
\item 
By Proposition \ref{Th:Prop_PBFTPK_1} and the second multiplicity condition, we have
\begin{equation}
\begin{split}
C_{1} &= \dfrac{\kappa_{a}\gamma}{\kappa_{a} - \kappa_{e}} y_{10} + x_{10}
\\
C_{2} &= \dfrac{\kappa_{a}\gamma}{\kappa_{a} - \kappa_{e}} y_{10} 
\end{split}    
\end{equation}
But, $y_{10} = d_{n}$ and $x_{10} = \overset{\scriptstyle{(n)}}{x_{1}}(t_{n-1}) = \overset{\scriptstyle{(n-1)}}{x_{2}}(t) = \mathbf{\overset{(n-1)}{Rem^{2}_{x}}}$\\
\item 
By the fourth multiplicity condition, we have
\begin{equation}
x(t) = x_{20}e^{-\kappa_{e}(t - s_{n})}    
\end{equation}
and Proposition \ref{Th:Prop_PBFTPK_2} implies that $x_{20} = \overset{\scriptstyle{(n-1)}}{x_{2}}(s_{n}) = \mathbf{\overset{(n)}{Rem^{1}_{x}}}$
\end{itemize}
% ----
\item
We now proceed to prove the formulas for the remains given in \eqref{Eq:FAT_IVP_REMAINS2} as follows:
\begin{itemize}
\item 
To prove the second formula for the remains, we only use the definition,
\begin{equation}
\begin{split}
\mathbf{\overset{(n)}{Rem^{1}_{x}}} &=
\overset{\scriptstyle{(n)}}{x_{2}}(t_{n})
= 
C_{3}e^{-\kappa_{e}(t_{n} - s_{n})}
\\
&= 
\mathbf{\overset{(n)}{Rem^{1}_{x}}}
\dfrac{e^{-\kappa_{e} t_{n}}}{e^{-\kappa_{e} s_{n}}}
\\
&= 
\mathbf{\overset{(n)}{Rem^{1}_{x}}}
\dfrac{\beta^{n}}{\beta^{n-1} B}
=
\mathbf{\overset{(n)}{Rem^{1}_{x}}}
\dfrac{\beta}{B}
\end{split}    
\end{equation}
\item 
To prove the formula for the first remainder, we proceed using mathematical induction:\\
\begin{enumerate}[-]
\item 
\textbf{Base case.} Verifying the formula for $n = 1$ is straightforward. It can be seen from Proposition \ref{Th:Prop_PBFTPK_1} by replacing the initial conditions.
\item 
\textbf{Induction hypothesis:}
Suppose the formula is valid for $n=k$
\begin{equation}
\mathbf{\overset{(k)}{Rem^{1}_{x}}} =
\dfrac{\kappa_{a}\cdot\gamma}{\kappa_{a} - \kappa_{e}}
(B - A)
\left[
\sum\limits_{i = 1}^{k}
\beta^{k-i}\,
d_{i}
\right]
\end{equation}
\item 
\textbf{Induction step:}
We must demonstrate that the formula is valid for $n=k+1$
\begin{equation}
\mathbf{\overset{(k+1)}{Rem^{1}_{x}}} =
\dfrac{\kappa_{a}\cdot\gamma}{\kappa_{a} - \kappa_{e}}
(B - A)
\left[
\sum\limits_{i = 1}^{k+1}
\beta^{(k+1)-i}\,
d_{i}
\right]
\end{equation}
Using the definition and Proposition \ref{Th:Prop_PBFTPK_2}, we have:
\begin{equation}
\begin{split}
\mathbf{\overset{(k+1)}{Rem^{1}_{x}}} &=
\overset{\scriptstyle{(k+1)}}{x_{2}}(s_{k+1})
\\
&=
\left[
\dfrac{\kappa_{a}\cdot\gamma}{\kappa_{a} - \kappa_{e}}
d_{k+1} + \overset{\scriptstyle{(k)}}{x_{2}}(t_{k})
\right]
e^{-\kappa_{e}(s_{k+1} - t_{k})}
-
\left[
\dfrac{\kappa_{a}\cdot\gamma}{\kappa_{a} - \kappa_{e}}
d_{k+1}
e^{-\kappa_{a}(s_{k+1} - t_{k})}
\right]
\\
&=
\left[
\dfrac{\kappa_{a}\cdot\gamma}{\kappa_{a} - \kappa_{e}}
d_{k+1} + \overset{\scriptstyle{(k)}}{x_{2}}(t_{k})
\right]
\dfrac{\beta^{k} B}{\beta^{k}}
-
\left[
\dfrac{\kappa_{a}\cdot\gamma}{\kappa_{a} - \kappa_{e}}
d_{k+1}
\right]
\dfrac{\alpha^{k} A}{\alpha^{k}}
\\
&=
\dfrac{\kappa_{a}\cdot\gamma}{\kappa_{a} - \kappa_{e}}
(B - A)d_{k+1}
+
\overset{\scriptstyle{(k)}}{x_{2}}(t_{n}) B
\\
&=
\dfrac{\kappa_{a}\cdot\gamma}{\kappa_{a} - \kappa_{e}}
(B - A)d_{k+1}
+
\mathbf{\overset{(k)}{Rem^{2}_{x}}}\,
B
\\
&=
\dfrac{\kappa_{a}\cdot\gamma}{\kappa_{a} - \kappa_{e}}
(B - A)d_{k+1}
+
\mathbf{\overset{(k)}{Rem^{1}_{x}}}\, B\,
\dfrac{\beta}{B}
\\
&=
\dfrac{\kappa_{a}\cdot\gamma}{\kappa_{a} - \kappa_{e}}
(B - A)d_{k+1}
+
\beta \, 
\dfrac{\kappa_{a}\cdot\gamma}{\kappa_{a} - \kappa_{e}}
(B - A)
\left[
\sum\limits_{i = 1}^{k}
\beta^{k-i}\,
d_{i}
\right]
\\
&=
\dfrac{\kappa_{a}\cdot\gamma}{\kappa_{a} - \kappa_{e}}
(B - A)
\left[
\sum\limits_{i = 1}^{k+1}
\beta^{(k+1)-i}\,
d_{i}
\right]
\end{split}
\end{equation}
\end{enumerate}
\end{itemize}
\end{enumerate}  
\end{proof}
\end{proof}

\end{document}